\documentclass[10pt]{amsart}
\usepackage{amsfonts}
\usepackage{amssymb}
\usepackage[latin1]{inputenc}
\usepackage{amsmath}
\usepackage{delarray}
\newtheorem{thm}{Theorem}[section]

\newtheorem{lemma}[thm]{Lemma}
\newtheorem{prop}[thm]{Proposition}
\theoremstyle{definition}

\newtheorem{defn}[thm]{Definition}
\newtheorem{remark}[thm]{Remark}
\numberwithin{equation}{section}
\newtheorem{example}[thm]{Example}

\newcommand{\fg}{\frak g}
\newcommand{\ft}{\frak t}

\newcommand{\cA}{\mathcal A}
\newcommand{\cB}{\mathcal B}

\newcommand{\cF}{\mathcal F}
\newcommand{\cH}{\mathcal H}
\newcommand{\cL}{\mathcal L}

\newcommand{\cR}{\mathcal R}
\newcommand{\cU}{\mathcal U}

\newcommand{\cX}{\mathcal X}

\newcommand{\bbR}{\mathbb R}
\newcommand{\bbT}{\mathbb T}
\newcommand{\bbC}{\mathbb C}
\newcommand{\bbN}{\mathbb N}
\newcommand{\bbK}{\mathbb K}
\newcommand{\bbS}{\mathbb S}
\newcommand{\bbZ}{\mathbb Z}

\newcommand{\ddim}{{\rm ddim\ }}
\newcommand{\rank}{{\rm rank\ }}

\newcommand{\QED}{\hfill $\square$\vspace{2mm}}
\newcommand{\Proof}{{\bf Proof}}

\begin{document}

\title{Torus Actions and Integrable Systems}

\author{Nguyen Tien Zung}
\address{Laboratoire Emile Picard, UMR 5580 CNRS, UFR MIG, Universit\'{e} Toulouse III}
\email{tienzung@picard.ups-tlse.fr}

\date{second version, July 2004}


\begin{abstract}
This is a survey on natural local torus actions which arise in integrable
dynamical systems, and their relations with other subjects, including:
reduced integrability, local normal forms, affine structures, monodromy,
global invariants, integrable surgery, convexity properties of momentum
maps, localization formulas, integrable PDEs.
\end{abstract}
\maketitle

\tableofcontents

\parskip=4pt

\section{Introduction}

To say that everything is a torus would be a great exaggeration, but to
say that {\it everything contains a torus} would not be too far from the
truth. According to ancient oriental philosophy, everything can be
described by (a combination of) five elemental aspects, or phases:
regular, transitive, expansive, chaotic, and contractive, and if we look
at these five phases as a whole then they also form cycles.

This survey paper is concerned with regular aspects of things.
Mathematically, they correspond to regular dynamics, or integrable
dynamical systems. For me, {\it an integrable system is a local torus
action}. The main dynamical property of a regular dynamical system is its
quasi-periodic behavior. Mathematically, it means that locally there is a
torus action which preserves the system. These torus actions exist near
compact regular orbits (Liouville's theorem). To a great extent, they
exist near singularities of integrable systems as well, and the main
topics of this paper are how to find them and what are the implications of
their existence.

I spent the last fifteen years looking for tori and this paper is mostly a
report on what I found -- probably very little for a fifteen year work. My
lame excuse is that I often had an empty stomach (empty pocket) to the
point of wanting to forget completely about the tori on more than one
occasions. This ``toric business'' began in 1989 when I was a 3rd year
undergraduate student under the direction of A.T. Fomenko: for my year-end
memoir I studied integrable perturbations of integrable Hamiltonian
systems of 2 degrees of freedom, and found out that there is a local
Hamiltonian $\bbT^1$-action which preserves an integrable 2D system near
each corank 1 nondegenerate singular hyperbolic level set
\cite{Zung-Bott1990}. This little discovery
is the starting point for more general existence results obtained
later.

The topics discussed in this paper can be seen from the table of contents.
They include: reduced integrability, Poincaré--Birkhoff normalization,
automorphism groups, partial action-angle variables, classification of
singularities, monodromy, characteristic classes, integrable surgery,
convexity of momentum maps, localization formulas, and so on. We consider
only classical dynamical systems in this paper. It turns out that many
local and semi-local results about the behavior of classical integrable
systems have their counterparts in quantum integrable systems, or at least
are useful for the study of quantum systems, see Vu Ngoc San
\cite{San-Semiclassic2004} and references therein.

The present paper deals mainly with finite-dimensional dynamical systems,
i.e. ordinary differential equations, though in the last section we will
briefly discuss the infinite dimensional case, i.e. integrable PDEs, where
there is a huge amount of literature but at the same time many basic
questions on topological aspects remain open.

{\bf Acknowledgements}. I would like to thank my Russian
colleagues Alexei Bolsinov, Anatoly Fomenko and Andrey Oshemkov for the
invitation to write this paper.

Many people helped me find the tori. I'm indebted to them all, and
especially to Michèle Audin, Alexei Bolsinov, Yves Colin de Verdière,
Richard Cushman, Jean-Paul Dufour, Hans Duistermaat, Anatoly Fomenko,
Lubomir Gavrilov, Thomas Kappeler, Pierre Molino, Tudor Ratiu, Jean-Claude
Sikorav, Vu Ngoc San, and Alan Weinstein.

A part of this paper was written during the author's visit to Bernoulli
Center, EPFL, Lausanne, and he would like to thank Tudor Ratiu and the
Bernoulli Center for the hospitality.

\section{Integrability, torus actions, and reduction}

\subsection{Integrability à la Liouville} \hfill

Probably the most well-known notion of integrability in dynamical systems
is the notion of integrability à la Liouville for Hamiltonian systems on
symplectic manifolds. Denote by $(M^{2n},\omega)$ a symplectic manifold of
dimension $2n$ with symplectic form $\omega$, and $H$ a function on
$M^{2n}$. Denote by $X_H$ the Hamiltonian vector field of $H$ on $M^{2n}$:
\begin{equation}
i_{X_H} \omega = -dH \ .
\end{equation}

\begin{defn}
A function $H$ (or the corresponding Hamiltonian vector field $X_H$) on a
$2n$-dimensional symplectic manifold $(M^{2n},\omega)$ is called {\bf
integrable à la Liouville}, or {\bf Liouville-integrable}, if it admits
$n$ functionally independent first integrals in involution. In other
words, there are $n$ functions $F_1 =H, F_2,\dots,F_n$ on $M^{2n}$ such
that $dF_1 \wedge \dots \wedge dF_n \neq 0$ almost everywhere and
$\{F_i,F_j\} = 0 \ \forall \ i,j$.
\end{defn}

In the above definition, $\{F_i,F_j\} := X_{F_i}(F_j)$ denotes the Poisson bracket
of $F_i$ and $F_j$ with respect to the symplectic form $\omega$. The map
\begin{equation} {\bf F}= (F_1,\dots,F_n): (M^{2n},\omega) \rightarrow
\bbK^n
\end{equation}
is called the {\bf momentum map} ($\bbK = \bbR$ or $\bbC$). The above
definition works in many categories: smooth, real analytic, holomorphic,
formal, etc.

The condition $X_H(F_i) = \{H,F_i\} = 0$
implies that the Hamiltonian vector field $X_H$ is tangent to the level
sets of $\bf F$. Let $N = {\bf F}^{-1}(c)$ be a regular connected
(component of a) level set of $\bf
F$. Then it is a Lagrangian submanifold of $M^{2n}$: the dimension of $N$ is
half the dimension of $M^{2n}$, and the restriction of
$\omega$ to $N$ is zero. So we can talk about
a (singular) Lagrangian foliation/fibration given by the momentum map.

A classical result attributed to Liouville \cite{Liouville-Torus1855} says
that, in the smooth case, if a connected level set $N$ is compact and does
not intersect with the boundary of $M^{2n}$, then it is diffeomorphic to a
standard torus $\bbT^n$, and the Hamiltonian system $X_H$ is
quasi-periodic on $N$: in other words, there is a periodic coordinate
system $(q_1,\dots,q_n)$ on $N$ with respect to which the restriction of
$X_H$ to $N$ has constant coefficients: $X_H = \sum \gamma_i \partial
/\partial q_i$, $\gamma_i$ being constants. For this reason, $N$ is called
a {\bf Liouville torus}.

The description of a Liouville-integrable Hamiltonian system near a
Liouville torus is given by the following theorem about the existence of
action-angle variables. This theorem is often called {Arnold-Liouville
theorem}, but it was essentially obtained by Henri Mineur in 1935
\cite{Mineur-AA1935,Mineur-AA1937}:

\begin{thm}[Liouville--Mineur--Arnold]
Let $N$ be a Liouville torus of a Liouville-integrable Hamiltonian
system with a given momentum map ${\bf F}: (M^{2n},\omega) \to
\bbR^n$. Then there is a neighborhood $\cU(N)$ of $N$ and a smooth
symplectomorphism
\begin{equation}
\Psi: (\cU(N),\omega) \to (D^n \times \bbT^n, \sum_1^n d\nu_i \wedge d\mu_i)
\end{equation}
($\nu_i$ - coordinates of $D^n$, $\mu_i \ (mod\ 1)$ - periodic
coordinates of
$\bbT^n$) such that $\bf F$ depends only on $I_i = \phi^\ast \nu_i$, i.e. $\bf F$
does not depend on $\phi_i = \phi^\ast \mu_i$.
\end{thm}
The variables $(I_i,\phi_i)$ in the above theorem are called {\bf
action-angle variables}. The map
\begin{equation}
(I_1,\dots,I_n) : (\cU(N),\omega) \to \bbR^n
\end{equation}
is the momentum map of a Hamiltonian torus $\bbT^n$-action on
$(\cU(n),\omega)$ which preserves $\bf F$. The existence of this
Hamiltonian torus action is essentially equivalent to
Liouville-Mineur-Arnold theorem: once the
action variables are found, angle variables can also be found easily by
fixing a Lagrangian section to the foliation by Liouville tori.
The quasi-periodicity of the system on $N$ also follows
immediately from the existence of this torus action.

The existence of action-angle variables is very important, both for the
theory of near-integrable systems (K.A.M. theory), and for the
quantization of integrable systems (Bohr--Sommerfeld rule). Actually,
Mineur was an astrophysicist, and Bohr--Sommerfeld quantization was his
motivation for finding action-angle variables.

Mineur \cite{Mineur-AA1937} also wrote down the following simple formula,
which we will call {\bf Mineur--Arnold formula}, for action functions:
\begin{equation}
\label{eqn:Mineur}
I_i (z) = \int_{\Gamma_i(z)} \beta
\end{equation}
where $z$ is a point in $\cU(N)$, $\beta$ is a primitive of the symplectic
form $\omega$, i.e. $d\beta = \omega$, and $\Gamma_i(z)$ is an 1-cycle on
the Liouville torus which contains $z$ (and which depends on $z$
continuously).

In the case of {\it algebraically integrable systems} (see e.g.
\cite{AdMoVa-Integrable2003}), where invariant tori can be identified with
(the real part of) Jacobian or Prym varieties of complex
curves ({\it spectral curves} of the system),
the integral in Mineur--Arnold formula corresponds to
Abelian integrals on complex curves,
as observed by Novikov and Veselov \cite{VeNo-Tori1984}.

It often happens that the above Mineur-Arnold formula (\ref{eqn:Mineur})
is still valid when the cycle $\Gamma_i$ lies on a singular level set, and
it leads to an action function near a singularity (see e.g.
\cite{Francoise-Period1990, Ito-AA1991} and Section \ref{section:local}).

\subsection{Generalized Liouville integrability} \hfill

In practice, one often deals with Hamiltonian systems which admit a
non-Abelian group of symmetries, or Hamiltonian systems on Poisson
(instead of symplectic) manifolds. A typical example is the Euler equation
on the dual of a Lie algebra. For such systems, Liouville integrability
needs to be replaced by a more general and convenient notion of
integrability, which nevertheless retains the main feature of Liouville
integrability, namely the existence of local torus actions.

Let $(M,\Pi)$ be a {\bf Poisson manifold}, with $\Pi$ being the Poisson
structure. It means that $\Pi$ is a 2-vector field on $M$ such that the
following binary operation on the space of functions on $M$, called the
{\bf Poisson bracket},
\begin{equation}
\{H,F\} = \langle dH \wedge dF, \Pi \rangle
\end{equation}
is a Lie bracket, i.e. it satisfies the Jacobi identity. A symplectic
manifold is also a Poisson manifold. Conversely, a Poisson manifold can be
seen as a singular foliation by symplectic manifolds, see e.g.
\cite{Weinstein-Review1998}.

Let $H$ be a function on a Poisson manifold $(M,\Pi)$, and $X_H$ the
corresponding Hamiltonian vector field: $X_H = dH \lrcorner \Pi$. Let
$\cF$ be a set of first integrals of $X_H$, i.e. each $F \in \cF$ is a
function on $M$ which is preserved by $X_H$ (equivalently, $\{F, H\} =
0$). Denote by $\ddim \cF$ the {\bf functional dimension} of $\cF$, i.e.
the maximal number of functions in $\cF$ which are functionally
independent almost everywhere.

We will associate to $\cF$ the space $\cX_\cF$ of Hamiltonian vector fields
$X_E$ such that $X_E(F) = 0$ for all $F \in \cF$ and $E$ is functionally dependent
of $\cF$ (i.e. the functional dimension of the union of $\cF$ with the function $E$
is the same as the functional dimension of $\cF$). Clearly, the vector fields in
$\cX_\cF$ commute pairwise and commute with $X_H$. Denote by $\ddim \cX_\cF$ the functional
dimension of $\cX_\cF$, i.e. the maximal number of vector fields in $X$
which are linearly independent at almost every point.
Note that we always have $\ddim \cF + \ddim \cX_\cF \leq m$,
because the vector fields in $\cX_\cF$ are tangent to the common level sets of the
functions in $\cF$.

The following definition is essentially due to
Nekhoroshev \cite{Nekhoroshev-Integrable1972} and Mischenko and Fomenko
\cite{MiFo-Liouville1978} :
\begin{defn}
\label{defn:HI} A Hamiltonian vector field $X_H$ on an $m$-dimensional
Poisson manifold $(M,\Pi)$ is called {\bf integrable in generalized
Liouville sense} if there is a set of first integrals $\cF$ such that
$\ddim \cF + \ddim \cX_\cF = m$.
\end{defn}

The above notion of integrability is also called {\bf noncommutative
integrability}, due to the fact that the functions in $\cF$ do not
Poisson-commute in general, and in many cases one may choose $\cF$ to be a
finite-dimensional non-commutative Lie algebra of functions under the
Poisson bracket. When the functions in $\cF$ Poisson-commute and the
Poisson structure is nondegenerate, we get back to the usual integrability
à la Liouville.

Denote $q = \ddim \cF, p = \ddim \cX_\cF$. Then we can find $p$
Hamiltonian vector fields $X_1 = X_{E_1},...,X_p = X_{E_p} \in \cX_\cF$ and $q$
functions $F_1,...,F_q \in \cF$ such that we have:
\begin{equation}
\begin{array}{c}
X_H(F_i) = 0,\ [X_H,X_i] = 0, \ [X_i,X_j] = 0, \ X_i(F_j) = 0 \ \forall \ i,j \ , \cr
X_1 \wedge \dots X_p \neq 0 \ {\rm and} \ dF_1 \wedge \dots dF_q \neq 0 \
{\rm almost \ everywhere} .
\end{array}
\end{equation}
The existence of such a $p$-tuple ${\bf X} = (X_1,...,X_p)$ of commuting
Hamiltonian vector fields and $q$-tuple ${\bf F} = (F_1,..., F_q)$ of
common first integrals with $p+q=m$ is equivalent to the integrability in
the generalized Liouville sense. When $p+q = m$, we will say that $H$ is
integrable with the aid of $({\bf X}, {\bf F})$, and by abuse of language,
we will also say that $({\bf X},{\bf F})$ is an {\bf integrable
Hamiltonian system} in generalized Liouville sense. The map
\begin{equation}
{\bf F} = (F_1,...,F_q): (M,\Pi) \to \bbK^q
\end{equation}
(where $\bbK = \bbR$ or $\bbC$) is called the {\bf generalized momentum
map}. The (regular) level sets of this map are called {\bf invariant
manifolds}: they are invariant with respect to $X_H$, $\bf X$ and $\bf F$.
They are of dimension $p$, lie on the symplectic leaves of $M$, and are
isotropic. When $p < {1 \over 2} \rank \Pi$, i.e. when the invariant
manifolds are isotropic but not Lagrangian, one also speaks of {\bf
degenerate integrability}, or {\bf superintegrability}, see e.g.
\cite{Evans-Superintegrability1983,
Nekhoroshev-Integrable1972,Reshetikhin-PoissonLie2003}.

\begin{defn}
With the above notations, a Hamiltonian system $X_H$, on a real Poisson
manifold $(M,\Pi)$, integrable with the aid of $({\bf X}, {\bf F})$, is
called {\bf proper} if the generalized momentum map ${\bf F} : M \to
\bbR^q$ is a proper map from $M$ to its image, and the image of the
singular set $\{x \in M , X_1 \wedge X_2 \wedge ... \wedge X_p (x) = 0\}$
of the commuting Hamiltonian vector fields under the momentum map ${\bf
F}: M \to \bbR^q$ is nowhere dense in $\bbR^q$.
\end{defn}

Under the properness condition, one gets a natural generalization of the
classical Liouville-Mineur-Arnold theorem
\cite{Nekhoroshev-Integrable1972,MiFo-Liouville1978}: outside the singular
region, the Poisson manifold $M$ is foliated by invariant isotropic
$p$-dimensional tori on which the flow of $X_H$ is quasi-periodic, and
there exist local action-angle coordinates. The action variables can still
be defined by Mineur-Arnold formula (\ref{eqn:Mineur}). There will be $p$
action and $p$ angle variables (so one will have to add $(q-p)$ variables
to get a full system of variables). In particular, near every isotropic
invariant torus there is a free Hamiltonian $\bbT^p$-action which
preserves the system.

\begin{example} A Hamiltonian $\bbT^p$-action on a Poisson manifold can be
seen as a proper integrable system -- the space of first integrals is the
space of $\bbT^p$-invariant functions, and in this case we have a {\it
global} $\bbT^p$-action which preserves the system. More generally, one
can associate to each Hamiltonian action of a compact Lie group $G$ on a
Poisson manifold a proper integrable system: $H$ is the composition of the
momentum map $\mu: M \rightarrow \fg^*$ with a generic $Ad^*$-invariant
function $h: \fg^* \rightarrow \bbR$, see Subsection
\ref{subsection:reduced1}.
\end{example}

\begin{remark}
There is a natural question: is an integrable Hamiltonian system in
generalized Liouville sense on a symplectic manifold also integrable à la
Liouville? One often expects the answer to be Yes. See e.g. Fomenko
\cite{Fomenko-Integrability1988} for a long discussion on this question,
and the related question about the existence of Liouville-integrable
systems on given symplectic manifolds.
\end{remark}

\begin{remark}
\label{remark:F_H}
Another natural question is the following. Let $\cF_H$
denote the space of all first integrals of $H$. Suppose that $H$ is
integrable in generalized Liouville sense. Is it true that $H$ is
integrable with the aid of $(\cF_H,\cX_{\cF_H})$, i.e. $\ddim \cF_H +
\ddim \cX_{\cF_H} = 0 ?$ We expect the answer to be yes for ``reasonable''
systems. It is easy to see that the answer is Yes
in the proper integrable case, under the additional assumption that the orbits
of $X_H$ are dense (i.e. its frequencies are incommensurable) on almost every
invariant torus (i.e. common level of a given set of first integrals $\cF$).
In this case $\cX_{\cF_H}$ consists of the Hamiltonian vector fields whose
flow is quasi-periodic on each invariant torus.
Another case where the answer is also Yes arises in the study of
local normal forms of analytic integrable vector fields, see Section
\ref{section:local}.
\end{remark}

\subsection{Non-Hamiltonian integrability} \hfill

There are many physical non-Hamiltonian (e.g. non-holonomic) systems,
which may naturally be called integrable in a non-Hamiltonian sense,
because their behavior is very similar to that of integrable Hamiltonian
systems, see e.g. \cite{BaCu-Nonholonomic1999,CuDu-Focus2001}. A simple
example is the Chinese top. (It is a spinning top whose lower part looks
like a hemisphere and whose upper part is heavy. When you spin it, it will
turn upside down after a while). The notion of non-Hamiltonian
integrability was probably first introduced by Bogoyavlenskij (see
\cite{Bogoyavlenskij-Integrability1998} and references therein), who calls
it {\bf broad integrability}, though other authors also arrived at it
independently, from different points of view, see e.g.
\cite{BaCu-Nonholonomic1999,Bogoyavlenskij-Integrability1998,CuDu-Focus2001,
Stolovitch-Singular2000,Zung-PD2002}.

\begin{defn}
A vector field $X$ on a  manifold $M$ is called {\bf integrable in
non-Hamiltonian sense} with the aid of $(\cF,\cX)$, where $\cF$ is a set
of functions on $M$ and
$\cX$ is a set of vector fields on $M$, if the following conditions are satisfied : \\
a)  $X(F) = 0 \ {\rm and} \ Y(F) = 0 \  \forall \ \ F \in \cF, Y \in \cX ,$ \\
b)  $[Y,X] = [Y,Z] = 0 \ \ \forall \ Y,Z \in \cX  ,$ \\
d) $\dim M = \ddim \cF + \ddim \cX .$ \\
In the real case, if, moreover, there is a $p$-tuple ${\bf X} =
(X_1,...,X_p)$ of vector fields in $\cX$ and a $q$-tuple ${\bf F} =
(F_1,...,F_q)$ of functionally independent functions in $\cF$, where $p =
\ddim \cX$ and $q = \ddim \cF$, such that the map ${\bf F} : M \to \bbR^q$
is a proper map from $M$ to its image, and for almost every level set of
this map the vector fields $X_1,...,X_p$ are linearly independent
everywhere on this level set, then we say that $X$ is {\bf proper
integrable} with the aid of $({\bf X},{\bf F})$, and by abuse of language
we will also say that $({\bf X},{\bf F})$ is a proper integrable
non-Hamiltonian system of bi-degree $(p,q)$ of freedom.
\end{defn}

So non-Hamiltonian integrability is almost the same as Hamiltonian
integrability, except for the fact that the vector fields
$X,X_1,\dots,X_p$ are not required to be Hamiltonian. It is not surprising
that Liouville's theorem holds for proper non-Hamiltonian integrable
systems as well: each regular invariant manifold (connected level set of
$\bf F$) is a $p$-dimensional torus on which the system is quasi-periodic,
and in a neighborhood of it there is a free $\bbT^p$-torus action which
preserves the system.

If a Hamiltonian system is (proper) integrable in the generalized
Liouville sense, then of course it is also (proper) integrable in the
non-Hamiltonian sense, though the inverse is not true: it may
happen that the invariant tori are not isotropic, see e.g.
\cite{Bogoyavlenskij-Integrability1998}.

\begin{remark}
\label{remark:F_X}
Remark \ref{remark:F_H} also applies to non-Hamiltonian systems:
For an integrable vector field $X$ on a manifold $M$,
denote by $\cF_X$ the set of all first integrals of $X$,
and by $\cX_X$ the set of all vector fields which commute with $X$
and preserve every function in $\cF$. Then a natural question is,
do we have the equality $\ddim \cF_X + \ddim \cX_X = \dim M$ ? The
answer is similar to the Hamiltonian case. In particular, if the system is proper
and the vector field $X$ is nonresonant (i.e. has a dense orbit) on almost
every invariant torus, then the answer is yes.
\end{remark}

\begin{remark} If $({\bf X}, {\bf F})$ is a non-Hamiltonian integrable
system of bidegree $(p,q)$ on a manifold $M$, then it can be lifted to a
Liouville-integrable system on the cotangent bundle $T^*M$. Denote by
$\pi: T^*M \rightarrow M$ the projection, then the corresponding momentum
map is $(H_1,\hdots,H_{p+q}): T^*M \rightarrow \bbR^{p+q}$, where $H_i =
\pi^* F_i$ ($i=1,\hdots,q$), and $H_{i+q} (\alpha) = \langle \alpha, X_i
(\pi(\alpha)) \rangle \ \forall\ \alpha \in T^*M$ ($i=1,\hdots,p)$.
\end{remark}

\subsection{Reduced integrability of Hamiltonian systems} \hfill
\label{subsection:reduced1}

In the literature, when people speak about integrability of a dynamical
system, they often actually mean its {\bf reduced integrability}, i.e.
integrability of the reduced (with respect to a natural symmetry group
action) system. For example, consider an integrable spinning top (e.g. the
Kovalevskaya top). Its configuration space is $SO(3)$, so it is naturally
a Hamiltonian system with 3 degrees of freedom. But it is often considered
as a 2-degree-of-freedom integrable system with a parameter, see e.g.
\cite{BoFo-Integrable1999}.

Curiously, to my knowledge, the natural question about the effect of
reduction on integrability has never been formally addressed in monographs
on dynamical systems. Recently we studied this question
\cite{Zung-Reduction2002}, and showed that, for a Hamiltonian system
invariant under a proper action of a Lie group, integrability is
essentially equivalent to reduced integrability.

It turns out that the most natural notion of integrability to use here is
not the Liouville integrability, but rather the integrability in
generalized Liouville sense. Also, since the category of manifolds is not
invariant under the operation of taking quotient with respect to a proper
group action, we have to replace manifolds by generalized manifolds: in
this paper, a {\bf generalized manifold} is a differentiable space which
is locally isomorphic to the quotient of a manifold by a compact group
action. Due to well-known results about functions invariant under compact
group actions, see e.g. \cite{Poenaru-Book1976}, one can talk about smooth
functions, vector fields, differential forms, etc. on generalized
manifolds, and the previous integrability definitions work for them as
well.

Let $(M,\Pi)$ be a Poisson generalized manifold, $G$ a Lie group
which acts properly on $M$, $H$ a function on $M$ which is invariant under
the action of $G$. Then the quotient space $M/G$ is again a Poisson
generalized manifold, see e.g. \cite{CuSj-Reduction1991}.
We will denote the projection of $\Pi,H,X_H$ on
$M/G$ by $\Pi/G,H/G, X_H/G$ respectively. Of course, $X_H/G$ is the
Hamiltonian vector field of $H/G$.

We will assume that the action of $G$ on $(M,\Pi)$ is Hamiltonian, with an
equivariant moment map $\mu : M \to {\frak g}^{\ast}$, where ${\frak g}$
denotes the Lie algebra of $G$, and that the following additional
condition is satisfied: Recall that the image $\mu(M)$ of $M$ under the
moment map $\mu : M \to {\frak g}^{\ast}$ is saturated by symplectic
leaves (i.e. coadjoint orbits) of ${\frak g}^{\ast}$. Denote by $s$ the
minimal codimension in ${\frak g}^{\ast}$ of a coadjoint orbit which lies
in $\pi(M)$. Then the additional condition is that there exist $s$
functions $f_1,...,f_s$ on ${\frak g}^{\ast}$, which are invariant on the
coadjoint orbits which lie in $\mu(M)$, and such that for almost every
point $x \in M$ we have $df_1 \wedge ... \wedge df_s (\mu(x)) \neq 0$. For
example, when $G$ is compact and $M$ is connected, then this condition is
satisfied automatically.

With the above notations and assumptions, we have :

\begin{thm}[\cite{Zung-Reduction2002}]
\label{thm:HI1} If the system $(M/G,X_H/G)$ is integrable in generalized
Liouville sense, then the system $(M,X_H)$ also is. Moreover, if $G$ is compact and $(M/G,
X_H/G)$
is proper, then $(M,X_H)$ also is.
\end{thm}

Since the preprint \cite{Zung-Reduction2002} will not be published as a
separate paper, let us include here a full proof of Theorem \ref{thm:HI1}.

\Proof. Denote by $\cF'$ a set of first integrals of $X_H/G$ on $M/G$
which provides the integrability of $X_H/G$, and by $\cX' = \cX_{\cF'}$
the corresponding space of commuting Hamiltonian vector fields on $M/G$.
We have $\dim M/G = p' + q'$ where $p' = \ddim \cX'$ and $q' = \ddim
\cF'$.

Recall that, by our assumptions, there exist $s$ functions $f_1,...,f_s$
on ${\frak g}^{\ast}$, which are functionally independent almost
everywhere in $\mu(M)$, and which are invariant on the coadjoint orbits
which lie in $\mu(M)$. Here $s$ is the minimal codimension in ${\frak
g}^{\ast}$ of the coadjoint orbits which lie in $\mu(M)$. We can complete
$(f_1,...,f_s)$ to a set of $d$ functions $f_1,...,f_s,f_{s+1},...,f_d$ on
${\frak g}^{\ast}$, where $d = \dim G = \dim \frak g$ denotes the
dimension of $\frak g$ , which are functionally independent almost
everywhere in $\mu(M)$.

Denote by $\overline{\cF}$ the pull-back of $\cF'$ under the projection
${\frak p} : M \to M/G $, and by $F_1,...,F_d$ the pull-back of
$f_1,...,f_d$ under the moment map $\mu : M \to {\frak g}^{\ast}$. Note
that, since $H$ is $G$-invariant, the functions $F_i$ are first integrals
of $X_H$. And of course, $\overline{\cF}$ is also a set of first integrals
of $X_H$. Denote by $\cF$ the union of $\overline{\cF}$ with
$(F_{s+1},...,F_d)$. (It is not necessary to include $F_1,...,F_s$ in this
union, because these functions are $G$-invariant and project to Casimir
functions on $M/G$, which implies that they are functionally dependent of
$\overline{\cF}$). We will show that $X_H$ is integrable with the aid of
$\cF$.

Notice that, by assumptions, the coadjoint orbits of ${\frak g}^{\ast}$
which lie in $\mu(M)$ are of generic dimension $d-s$, and the functions
$f_{s+1},...,f_d$ may be viewed as a coordinate system on a symplectic
leaf of $\mu(M)$ at a generic point. In particular, we have
$$\langle df_{s+1} \wedge ... \wedge df_d, X_{f_{s+1}} \wedge ...
X_{f_d} \rangle \neq 0,$$ which implies, by equivariance :
$$\langle dF_{s+1} \wedge ... \wedge dF_d, X_{F_{s+1}} \wedge ... X_{F_d}
\rangle \neq 0.$$

Since the vector fields $X_{F_{s+1}},...,X_{F_d}$ are tangent to the orbits of $G$
on $M$, and the functions in $\overline{\cF}$ are invariant on the orbits of $G$, it
implies that the set $(F_{s+1},...,F_{d})$ is ``totally'' functionally independent
of $\overline{\cF}$. In particular, we have :

\begin{equation}
\ddim \cF = \ddim \cF' + \ddim (F_{s+1},...,F_{d}) = q' + d - s ,
\end{equation}
where $q' = \ddim \cF'$. On the other hand, we have

$$\dim M = \dim M/G + (d-k) = p' + q' + d - k ,$$
where $p' = \ddim \cX_{\cF'}$, and $k$ is the dimension of a minimal isotropic group
of the action of $G$ on $M$. Thus, in order to show the integrability condition
$$\dim M = \ddim \cF + \ddim \cX_{\cF} ,$$ it remains to show that

\begin{equation}
\label{eqn:ddimcX}
\ddim \cX_{\cF} = \ddim \cX_{\cF'} + (s-k) .
\end{equation}

Consider the vector fields $Y_1 = X_{F_1},...,Y_d = X_{F_d}$ on $M$. They span the
tangent space to the orbit of $G$ on $M$ at a generic point. The dimension of such a
generic tangent space is $d-k$. It implies that, among the first $s$ vector fields,
there are at least $s-k$ vector fields which are linearly independent at a generic
points : we may assume that $Y_1 \wedge ... \wedge Y_{s-k} \neq 0 .$

Let $X_{h_1},...,X_{h_{p'}}$ be $p'$ linearly independent (at a generic point)
vector fields which belong to $\cX_{\cF'}$, where $p' = \ddim \cX_{\cF'}$. Then we
have

$$
X_{{\frak p}^{\ast}(h_1)},...,X_{{\frak p}^{\ast}(h_{p'})},Y_1,...,Y_{s-k} \in
\cX_{\cF},
$$
and these $p' + s - k$ vector fields are linearly independent at a generic point.
(Recall that, at each point $x \in M$, the vectors $Y_1(x),...,Y_{s-k}(x)$ are
tangent to the orbit of $G$ which contains $x$, while the linear space spanned by
$X_{{\frak p}^{\ast}(h_1)},...,X_{{\frak p}^{\ast}(h_{p'})}$ contains no tangent
direction to this orbit).

Thus we have $\ddim \cX_{\cF} \geq p' + s - k$, which means that
$\ddim \cX_{\cF} = p' + s - k$ (because, as discussed earlier, we
always have $\ddim \cF + \ddim \cX_\cF \leq \dim M$). We have
proved that if $(M/G,X_H/G)$ is integrable in generalized
Liouville sense then $(M,X_H)$ also is.

Now suppose that $G$ is compact and $(M/G,X_H/G)$ is proper:
there are $q'$ functionally independent functions $g_1,...,g_{q'} \in \cF'$ such
that $(g_1,...,g_{q'}): M/G \to \bbR^{q'}$ is a proper map from $M/G$ to its image,
and $p'$ Hamiltonian vector fields $X_{h_1},...,X_{h_{p'}}$ in $\cX'$ such that on a
generic common level set of $(g_1,...,g_{q'})$ we have that $X_{h_1} \wedge ...
\wedge X_{h_{p'}}$ does not vanish anywhere. Then it is straightforward that
$${\frak p}^{\ast}(g_1),...,{\frak p}^{\ast}(g_{q'}),F_{s+1},...,F_d \in \cF$$ and
the map $$({\frak p}^{\ast}(g_1),...,{\frak
p}^{\ast}(g_{q'}),F_{s+1},...,F_d) : M \to \bbR^{q'+d-s}$$ is a proper map
from $M$ to its image. More importantly, on a generic level set of this
map we have that the $(q'+s-k)$-vector $X_{{\frak p}^{\ast}(h_1)} \wedge
... \wedge X_{{\frak p}^{\ast}(h_{p'})} \wedge Y_1 \wedge ... \wedge
Y_{s-k}$ does not vanish anywhere. To prove this last fact, notice that
$X_{{\frak p}^{\ast}(h_1)} \wedge ... \wedge X_{{\frak p}^{\ast}(h_{p'})}
\wedge Y_1 \wedge ... \wedge Y_{s-k} (x) \neq 0$ for a point $x \in M$ if
and only if $X_{{\frak p}^{\ast}(h_1)} \wedge ... \wedge X_{{\frak
p}^{\ast}(h_{p'})} (x) \neq 0$ and $Y_1 \wedge ... \wedge Y_{s-k} (x) \neq
0$ (one of these two multi-vectors is transversal to the $G$-orbit of $x$
while the other one ``lies on it''), and that these inequalities are $G
\times \bbR^{p'}$-invariant properties, where the action of $\bbR^{p'}$ is
generated by $X_{{\frak p}^{\ast}(h_1)},...,X_{{\frak p}^{\ast}(h_{p'})}$.
\QED

\begin{remark}
Similar results to Theorem \ref{thm:HI1} have been obtained independently
by Bolsinov and Jovanovic
\cite{BoJo-Integrable2003,Jovanovic-Geodesic2002}, who used them to
construct new examples of integrable geodesic flows, e.g. on biquotients
of compact Lie groups.
\end{remark}

\begin{example} The simplest example which shows an evident
relationship between reduction and integrability is the classical Euler
top : it can be written as a Hamiltonian system on $T^{\ast}SO(3)$,
invariant under a natural Hamiltonian action of $SO(3)$, is integrable
with the aid of a set of four first integrals, and has 2-dimensional
isotropic invariant tori. The geodesic flow of a bi-invariant metric on a
compact Lie group is also properly integrable : in fact, the corresponding
reduced system is trivial (identically zero). More generally, let $H = h
\circ \mu$ be a collective Hamiltonian in the sense of
Guillemin--Sternberg (see e.g. \cite{GuSt-Collective1983}), where $\mu: M
\to \mathfrak{g}^*$ is the momentum map of a Hamiltonian compact group
action, and $h$ is a function on $\mathfrak{g}^*$. If $h$ is a Casimir
function on $\mathfrak{g}^*$, then $H$ is integrable because its reduction
will be a trivial Hamiltonian system.
\end{example}

\begin{remark} Recall from Equation (\ref{eqn:ddimcX}) that $\ddim
\cX_{\cF} - \ddim \cX_{\cF'} = s-k$, where $k$ is the dimension of a
generic isotropic group of the $G$-action on $M$, and $s$ is the (minimal)
corank in ${\frak g}^{\ast}$ of a coadjoint orbit which lies in $\mu(M)$.
On the other hand, the difference between the rank of the Poisson
structure on $M$ and the reduced Poisson structure on $M/G$ can be
calculated as follows :
\begin{equation}
\rank \Pi - \rank \Pi/G = (d-k) + (s-k)
\end{equation}
Here $(d-k)$ is the difference between $\dim M$ and $\dim M/G$, and $(s-k)$ is the
difference between the corank of $\Pi/G$ in $M/G$ and the corank of $\Pi$ in $M$. It
follows that
\begin{equation}
\rank \Pi - 2 \ddim \cX_{\cF} = \rank \Pi/G - 2 \ddim \cX_{\cF'} + (d-s)
\end{equation}
In particular, if $d-s > 0$ (typical situation when $G$ is non-Abelian), then we
always have $\rank \Pi - 2 \ddim \cX_{\cF} > 0$ (because we always have $\rank \Pi/G
- 2 \ddim \cX_{\cF'} \geq 0$ due to integrability), i.e. the original system is
always super-integrable with the aid of $\cF$. When $G$ is Abelian (implying $d=s$),
and the reduced system is Liouville-integrable with the aid of $\cF'$ (i.e. $\rank
\Pi/G = 2 \ddim \cX_{\cF'}$), then the original system is also Liouville-integrable
with the aid of $\cF$.
\end{remark}

\begin{remark} Following Mischenko-Fomenko \cite{MiFo-Liouville1978},
we will say that a
hamiltonian system $(M,\Pi,X_H)$ is {\bf non-commutatively integrable in
the restricted sense} with the aid of $\cF$ , if $\cF$ is a
finite-dimensional Lie algebra under the Poisson bracket and $(M,\Pi,X_H)$
is integrable with the aid of $\cF$. In other words, we have an
equivariant moment maps $(M,\Pi) \to {\frak f}^{\ast}$, where ${\frak f}$
is some finite-dimensional Lie algebra, and if we denote by $f_1,...,f_n$
the components of this moment map, then they are first integrals of $X_H$,
and $X_H$ is integrable with the aid of this set of first integrals.
Theorem \ref{thm:HI1} remains true, and its proof remains the same if not
easier, if we replace Hamiltonian integrability by non-commutative
integrability in the restricted sense. Indeed, if $M \to {\frak g}^{\ast}$
is the equivariant moment map of the symmetry group $G$, and if $M/G \to
{\frak h}^{\ast}$ is an equivariant moment map which provides
non-commutative integrability in the restricted sense on $M/G$, then the
map $M \to {\frak h}^{\ast}$ (which is the composition $M \to M/G
\to{\frak h}^{\ast}$) is an equivariant moment map which commutes with $M
\to {\frak g}^{\ast}$, and the direct sum of this two maps, $M \to {\frak
f}^{\ast}$ where $\frak f = \frak g \bigoplus \frak h$, will provide
non-commutative integrability in the restricted sense on $M$.
\end{remark}

Theorem \ref{thm:HI1} has the following inverse
(see Remark \ref{remark:F_H}):

\begin{thm}
\label{thm:HI2} If $G$ is compact, and if the Hamiltonian system $(M,X_H)$ is
integrable with the aid of $\cF_H$ (the set of all first integrals of $H$)
in the sense that $\ddim \cF_H \ +\ \ddim \cX_{\cF_H} = \dim M$,
then the reduced Hamiltonian system $(M/G,X_H)$ is also integrable.
Moreover, if $(M,X_H)$ is proper then $(M/G,X_H)$ also is.
\end{thm}

\Proof. By assumptions, we have $\dim M = p + q$, where $q = \ddim \cF_H$
and $p = \ddim \cX_{\cF_H}$, and we can find $p$ first integrals
$H_1,...,H_p$ of $H$ such that $X_{H_1},...,X_{H_p}$ are linearly
independent (at a generic point) and belong to $\cX_{\cF_H}$. In
particular, we have $X_{H_i}(F) = 0$ for any  $F \in \cF$ and $1 \leq i
\leq p$.

An important observation is that the functions $H_1,...,H_p$ are $G$-invariant. In
deed, if we denote by $F_1,...,F_d$ the components of the equivariant moment map
$\pi : M \to {\frak g}^{\ast}$ (via an identification of ${\frak g}^{\ast}$ with
$\bbR^d$), then since $H$ is $G$-invariant we have $\{H,F_j\}=0$, i.e. $F_j \in
\cF_H$, which implies that $\{F_j,H_i\} = 0 \ \ \forall 1 \leq i \leq d, \ \ 1 \leq
j \leq p$, which means that $H_i$ are $G$-invariant.

The Hamiltonian vector fields $X_{H_i}/G$ belong to
$\cX_{\cF_{H/G}}$ : Indeed, if $f \in \cF_{H/G}$ then ${\frak p}^{\ast}(f)$ is a first
integral of $H$, implying $\{H_i,{\frak p}^{\ast}(f)\} = 0$, or $\{H_i/G,f\}=0$, where
$\frak p$ denotes the projection $M \to M/G$.

To prove the integrability of $X_H/G$, it is sufficient to show that
\begin{equation}
\label{eqn:Xh} \dim M/G \leq \ddim \cF_{H/G} + \ddim (X_{H_1}/G,...,X_{H_q}/G)
\end{equation}
But we denote by $r$ the generic dimension of the intersection of a common level set
of $p$ independent first integrals of $X_H$ with an orbit of $G$ in $M$, then one
can check that
$$
p - \ddim (X_{H_1}/G,...,X_{H_q}/G) = \ddim \cX_{\cF_H} - \ddim (X_{H_1}/G,...,
X_{H_q}/G) = r
$$
and
$$ q - \ddim \cF_{H/G} = \ddim \cF_H - \ddim \cF_{H/G} \leq (d-k) - r
$$
where $(d-k)$ is the dimension of a generic orbit of $G$ in $M$. To prove the last
inequality, notice that functions in $\cF_{H/G}$ can be obtained from functions in
$\cF_H$ by averaging with respect to the $G$-action. Also, $G$ acts on the
(separated) space of common level sets of the functions in $\cF_H$, and isotropic
groups of this $G$-action are of (generic) codimension $(d-k) - r$.

The above two formulas, together with $p+q = \dim M = \dim M/G + (d-k)$,
implies Inequality (\ref{eqn:Xh}) (it is in fact an equality). The proper
case is straightforward. \QED

\subsection{Non-Hamiltonian reduced integrability}  \hfill

One of the main differences between the non-Hamiltonian case and the Hamiltonian
case is that reduced non-Hamiltonian integrability does not imply integrability. In
fact, in the Hamiltonian case, we can lift Hamiltonian vector fields from $M/G$ to
$M$ via the lifting of corresponding functions. In the non-Hamiltonian case, no such
canonical lifting exists, therefore commuting vector fields on $M/G$ do not provide
commuting vector fields on $M$. For example, consider a vector field of the type $X
= a_1 \partial/\partial x_1 + a_2 \partial/\partial x_2 + b(x_1,x_2)
\partial/\partial x_3$ on the standard torus $\bbT^3$ with periodic coordinates
$(x_1,x_2,x_3)$, where $a_1$ and $a_2$ are two incommensurable real numbers
($a_1/a_2 \notin {\mathbb Q}$), and $b(x_1,x_2)$ is a smooth function of two
variables. Then clearly $X$ is invariant under the ${\mathbb T}^1$-action generated
by $\partial/\partial x_3$, and the reduced system is integrable. On the other hand,
for $X$ to be integrable, we must be able to find a function $c(x_1,x_2)$ such that
$[X, \partial/\partial x_1 + c(x_1,x_2)\partial/\partial x_3] = 0 $. This last
equation does not always have a solution (it is a small divisor problem, and depends
on $a_1/a_2$ and the behavior of the coefficients of $b(x_1,x_2)$ in its Fourier
expansion), i.e. there are choices of $a_1,a_2,b(x_1,x_2)$ for which the vector
field $X$ is not integrable.

However, non-Hamiltonian integrability still implies reduced integrability.
Recall from Remark \ref{remark:F_X} that if a vector field $X$ on a
(generalized) manifold $M$ is integrable, then under mild additional
conditions we have $\ddim \cX_X + \ddim \cF_X = \dim M$, where
$\cF_X$ is the set of all first
integrals of $X$, and $\cX_X$ is the set of all vector fields which commute with $X$
and preserve every function in $\cF$.

\begin{thm}
\label{thm:nHI} Let $X$ be a smooth non-Hamiltonian proper integrable system on
a manifold $M$ with the aid of $(\cF_X,\cX_X)$, i.e. $\ddim
\cX_X + \ddim \cF_X = \dim M$, and $G$ be a compact Lie group
acting on $M$ which preserves $X$. Then the reduced system on $M/G$ is also proper
integrable.
\end{thm}

\Proof. Let $\cX_X^G$ denote the set of vector fields which belong to
$\cX_X$ and which are invariant under the action of $G$. Note that the
elements of $\cX_X^G$ can be obtained from the elements of $\cX_X$ by
averaging with respect to the $G$-action.

A key ingredient of the proof is the fact $\ddim \cX_X^G = \ddim \cX_X$ (To see this
fact, notice that near each regular invariant torus of the system there is an
effective torus action (of the same dimension) which preserves the system, and this
torus action must necessarily commute with the action of $G$. The generators of this
torus action are linearly independent vector fields which belong to  $\cX_X^G$ - in
fact, they are defined locally near the union of $G$-orbits which by an invariant
torus, but then we can extend them to global vector fields which lie in $\cX_X^G$)

Therefore, we can project the pairwise commuting vector fields in
$\cX_X^G$ from $M$ to $M/G$ to get pairwise commuting vector fields on
$M/G$. To get the first integrals for the reduced system, we can also take
the first integrals of $X$ on $M$ and average them with respect to the
$G$-action to make them $G$-invariant.  The rest of the proof of Theorem
\ref{thm:nHI} is similar to that of Theorem \ref{thm:HI2}.
\QED

\section{Torus actions and local normal forms}
\label{section:local}

\subsection{Toric characterization of Poincaré-Birkhoff normal form} \hfill
\label{subsection:action}

It is a simple well-known fact that every vector field near an equilibrium
point admits a formal Poincaré-Birkhoff normal form (Birkhoff in the
Hamiltonian case, and Poincaré-Dulac in the non-Hamiltonian case). What is
also very simple but much less well-known is the fact that these normal
forms are governed by torus actions. We will explain this fact here,
following \cite{Zung-Birkhoff2001,Zung-PD2002}.

Let $X$ be a given analytic vector field in a neighborhood of $0$ in
$\bbK^m$, where $\bbK = \bbR$ or $\bbC$, with $X(0) = 0$. When $\bbK = \bbR$,
we may also view $X$ as a holomorphic (i.e. complex analytic) vector field by
complexifying it. Denote by
\begin{equation}
X = X^{(1)} + X^{(2)} + X^{(3)} + ...
\end{equation}
the Taylor expansion of $X$ in some local system of coordinates, where $X^{(k)}$ is
a homogeneous vector field of degree $k$ for each $k \geq 1$.

In the Hamiltonian case, on a symplectic manifold, $X = X_H$, $m= 2n$,
$\bbK^{2n}$ has a standard symplectic structure, and $X^{(j)} =
X_{H^{(j+1)}}$.

The algebra of linear
vector fields on ${\mathbb K}^{m}$, under the standard Lie bracket, is nothing but
the reductive algebra $gl(m, \bbK) = sl(m,\bbK) \oplus \bbK$. In particular, we have
\begin{equation}
X^{(1)} = X^s + X^{nil},
\end{equation}
where $X^s$ (resp., $X^{nil}$) denotes the semi-simple (resp., nilpotent) part of
$X^{(1)}$. There is a complex linear system of coordinates $(x_j)$ in ${\mathbb
C}^{m}$ which puts $X^s$ into diagonal form:
\begin{equation}
X^s = \sum_{j=1}^m \gamma_j x_j \partial / \partial x_j ,
\end{equation}
where $\gamma_j$ are complex coefficients, called {\bf eigenvalues} of $X$
(or $X^{(1)}$) at $0$.

In the Hamiltonian case, $X^{(1)} \in sp(2n,\bbK)$ which is a simple Lie
algebra, and we also have the decomposition $X^{(1)} = X^s + X^{nil}$,
which corresponds to the decomposition
\begin{equation}
H^{(2)} = H^s + H^{nil}
\end{equation}
There is a complex canonical linear system of coordinates $(x_j,y_j)$ in ${\mathbb
C}^{2n}$ in which $H^s$ has diagonal form:
\begin{equation}
H^s = \sum_{j=1}^n \lambda_j x_j y_j ,
\end{equation}
where $\lambda_j$ are complex coefficients, called {\bf frequencies} of
$H$ (or $H^{(2)}$) at $0$.

For each natural number $k \geq 1$, the vector field $X^s$ acts linearly on the
space of homogeneous vector fields of degree $k$ by the Lie bracket, and the
monomial vector fields are the eigenvectors of this action:
\begin{equation}
[\sum_{j=1}^m \gamma_j x_j \partial / \partial x_j , x_1^{b_1}x_2^{b_2}...x_n^{b_n}
\partial / \partial x_l] = (\sum_{j=1}^n b_j\gamma_j - \gamma_l) x_1^{b_1}x_2^{b_2}...x_n^{b_n}
\partial / \partial x_l .
\end{equation}

When an equality of the type
\begin{equation}
\sum_{j=1}^m b_j\gamma_j - \gamma_l = 0
\end{equation}
holds for some nonnegative integer $m$-tuple $(b_j)$ with $\sum b_j \geq 2$, we will
say that the monomial vector field $x_1^{b_1}x_2^{b_2}...x_m^{b_m}
\partial / \partial x_l$ is a {\bf resonant term},
and that the $m$-tuple $(b_1,...,b_l - 1,..., b_l)$
is a resonance relation for the eigenvalues $(\gamma_i)$. More precisely,
a {\bf resonance relation} for  the $n$-tuple of eigenvalues $(\gamma_j)$
of a vector field $X$ is an $m$-tuple $(c_j)$ of integers satisfying the
relation $\sum c_j \gamma_j = 0,$ such that $c_j \geq -1, \sum c_j \geq
1,$ and at most one of the $c_j$ may be negative.

In the Hamiltonian case, $H^s$ acts linearly on the space of functions by
the Poisson bracket. Resonant terms (i.e. generators of the kernel of this
action) are monomials $\prod x_j^{a_j}y_j^{b_j}$ which satisfy the
following resonance relation, with $c_j = a_j - b_j$:
\begin{equation}
\sum_{j=1}^m c_j\lambda_j  = 0
\end{equation}

Denote by ${\mathcal R}$  the subset of ${\mathbb Z}^m$
(or sublattice of ${\mathbb Z}^n$ in the Hamiltonian case)
consisting of all resonance relations $(c_j)$ for a
given vector field $X$. The number
\begin{equation}
r = \dim_\bbZ (\cR \otimes \bbZ)
\end{equation}
is called the {\bf degree of resonance} of $X$. Of course, the degree of
resonance depends only on the eigenvalues of the linear part of $X$, and
does not depend on the choice of local coordinates. If $r=0$ then we say
that the system is {\bf nonresonant} at 0.

The vector field $X$ is said to be in {\bf Poincaré-Birkhoff normal form}
if it commutes with the semisimple part of its linear part (see e.g.
\cite{Bruno-Local1989,Roussarie-Asterisque1975}):
\begin{equation}
[X,X^s] = 0.
\end{equation}
In the Hamiltonian case, the above equation can also be written as
\begin{equation}
\{H,H^s\} = 0.
\end{equation}

The above equations mean that if $X$ is in normal form then its nonlinear
terms are resonant.  A transformation of coordinates (which is symplectic
in the Hamiltonian case) which puts $X$ in Poincaré-Birkhoff normal form
is called a {\bf Poincaré-Birkhoff normalization}. It is a classical
result of Poincaré, Dulac, and Birkhoff that any analytic vector field
which vanishes at 0 admits a {\it formal} Poincaré-Birkhoff normalization
(which does not converge in general).

Denote by ${\mathcal Q} \subset {\mathbb Z}^m$ the integral sublattice of ${\mathbb
Z}^m$ consisting of $m$-dimensional vectors $(\rho_j) \in {\mathbb Z}^m$ which
satisfy the following properties :
\begin{equation}
\label{eqn:Q} \sum_{j=1}^m \rho_j c_j = 0 \  \forall \ (c_j) \in {\mathcal R} \ , \
{\rm and} \ \ \rho_j = \rho_k \ \ {\rm if} \ \ \gamma_j = \gamma_k \
\end{equation}
(where $\mathcal R$ is the set of resonance relations as before). In the
Hamiltonian case, $\mathcal Q$ is defined by
\begin{equation}
\label{eqn:Q2} \sum_{j=1}^n \rho_j c_j = 0 \  \forall \ (c_j) \in {\mathcal R} \
.
\end{equation}

We will call the number
\begin{equation}
\label{eqn:t} d = \dim_\bbZ \mathcal Q
\end{equation}
the {\bf toric degree} of $X$ at $0$. Of course, this number depends only
on the eigenvalues of the linear part of $X$, and we have the following
(in)equality : $r + d = n$ in the Hamiltonian case (where $r$ is the
degree of resonance), and $r + d \leq m$ in the non-Hamiltonian case.

Let $(\rho^{1}_j),...,(\rho^d_j)$ be a basis of $\mathcal Q$. For each $k=1,...,
d$
define the following diagonal linear vector field $Z_k$ :

\begin{equation}
\label{eqn:Z} Z_k = \sum_{j=1}^m \rho^k_j x_j \partial / \partial x_j
\end{equation}
in the non-Hamiltonian case, and $Z_k = X_{F^k}$ where
\begin{equation}
\label{eqn:F} F^k = \sum_{j=1}^n \rho^k_j x_jy_j
\end{equation}
in the Hamiltonian case.

The vector fields $Z_1,...,Z_r$ have the following remarkable properties :

a) They commute pairwise and commute with $X^s$ and $X^{nil}$, and they are linearly
independent almost everywhere.

b) $iZ_j$ is a periodic vector field of period $2\pi$ for each $j \leq r$ (here $i =
\sqrt{-1}$). What does it mean is that if we write $iZ_j = \Re (iZ_j) + i \Im (iZ_j)
$, then $\Re (iZ_j)$ is a periodic real vector field in ${\mathbb C}^n = {\mathbb
R}^{2n}$ which preserves the complex structure.

c) Together, $iZ_1,..., iZ_r$ generate an effective linear ${\mathbb T}^r$-action in
${\mathbb C}^n$ (which preserves the symplectic structure in the Hamiltonian case),
which preserves $X^s$ and $X^{nil}$.

A simple calculation shows that $X$ is in Poincaré-Birkhoff normal form,
i.e. $[X,X^s] = 0$, if and only if we have
\begin{equation}
[X,Z_k] = 0 \ \ \ \forall \  k=1,...,r.
\end{equation}

The above commutation relations mean that if $X$ is in normal form, then
it is preserved by the effective $r$-dimensional torus action generated by
$iZ_1,...,iZ_r$. Conversely, if there is a torus action which preserves $
X$, then because the torus is a compact group we can linearize this torus
action (using Bochner's linearization theorem
\cite{Bochner-Linearization1945} in the non-Hamiltonian case, and the
equivariant Darboux theorem in the Hamiltonian case, see e.g.
\cite{CoDaMo-Moment1988,GuSt-Convexity1982}), leading to a normalization
of $X$. In other words, we have:

\begin{thm}[\cite{Zung-Birkhoff2001,Zung-PD2002}]
\label{thm:toricPB} A holomorphic (Hamiltonian) vector field $X$ in a neighborhood of $0$ in
${\mathbb C}^m$ (or $\bbC^{2n}$ with a standard symplectic form)
admits a locally holomorphic Poincaré-Birkhoff normalization if and only if it
is preserved by an effective holomorphic (Hamiltonian)
action of a real torus of dimension $t$,
where $t$ is the toric degree of $X^{(1)}$ as defined in (\ref{eqn:t}), in a
neighborhood of $0$ in ${\mathbb C}^m$ (or $\bbC^{2n}$),
which has $0$ as a fixed point and whose
linear part at $0$ has appropriate weights (given by the lattice $\mathcal Q$
defined in (\ref{eqn:Q},\ref{eqn:Q2}),
which depends only on the linear part $X^{(1)}$ of $X$).
\end{thm}

The above theorem is true in the formal category as well. But of
course, any vector field admits a formal Poincaré-Birkhoff normalization,
and a formal torus action.

\subsection{Some simple consequences and generalizations}    \hfill

Theorem \ref{thm:toricPB} has many important implications. One of them is:

\begin{prop}[\cite{Zung-Birkhoff2001,Zung-PD2002}]
\label{prop:realcomplex} A real analytic vector field $X$ (Hamiltonian
or non-Hamiltonian) in the neighborhood of an equilibrium point
admits a local real analytic Poincaré-Birkhoff normalization if and only
if it admits a local holomorphic Poincaré-Birkhoff
normalization when considered as a holomorphic vector field.
\end{prop}

The proof of the above proposition (see \cite{Zung-Birkhoff2001})
is based on the fact that the complex conjugation induces an involution on the
torus action which governs the
Poincaré-Birkhoff normalization.

If a dynamical system near an equilibrium point is invariant with respect
to a compact group action which fixes the equilibrium point, then this
compact group action commutes with the (formal) torus action of the
Poincaré-Birkhoff normalization. Together, they form a bigger compact
group action, whose linearization leads to a simultaneous
Poincaré-Birkhoff normalization and linearization of the compact symmetry
group, i.e. we can perform the Poincaré-Birkhoff normalization in an
invariant way. This is a known result in dynamical systems, see e.g.
\cite{Zhitomirskii-NF1994}, but the toric point of view gives a new simple
proof of it. The case of equivariant vector fields is similar. For
example, one can speak about Poincaré-Dulac normal forms for
time-reversible vector fields, see e.g. \cite{LaRo-Reversal1998}.

Another situation where one can use the toric characterization is the case
of isochore (i.e. volume preserving) vector fields. In this case,
naturally, the normalization transformation is required to be volume-preserving.
Both Theorem \ref{thm:toricPB} and Proposition \ref{prop:realcomplex} remain
valid in this case.

One can probably use the toric point of view to study normal forms of
Hamiltonian vector field on {\it Poisson} manifolds as well. For example,
let ${\mathfrak g}^*$ be the dual of a semi-simple Lie algebra, equipped
with the standard linear Poisson structure, and let $H: {\mathfrak g}^*
\to \bbK$ be a regular function near the origin 0 of ${\mathfrak g}^*$.
The corresponding Hamiltonian vector field $X_H$ will vanish at $0$,
because the Poisson structure itself vanishes at 0. Applying
Poincaré-Birkhoff normalization techniques, we can kill the ``nonresonant
terms'' in $H$ (with respect to the linear part of $H$, or $dH(0)$). The
normalized Hamiltonian will be invariant under the coadjoint action of a
subtorus of a Cartan torus of the (complexified) Lie group of $\mathfrak
g$. In the ``nonresonant'' case, we have a Cartan torus action which
preserves the system.

\subsection{Convergent normalization for integrable systems}  \hfill
\label{subsection:convergentintegrable}

Though every vector field near an equilibrium admits a formal
Poincaré-Birkhoff normalization, the problem of finding a convergent (i.e.
locally real analytic or holomorphic) normalization is much more
difficult. The usual step by step killing of non-resonant terms leads to
an infinite product of coordinate transformations, which may diverge in
general, due to the presence of small divisors. Positive results about the
convergence of this process are due to Poincaré, Siegel, Bruno and others
mathematicians, under Diophantine conditions on the eigenvalues of the
linear part of the system, see e.g.
\cite{Bruno-Local1989,Roussarie-Asterisque1975}.

However, when the vector field is analytically integrable (i.e. it is an
real or complex analytic vector field, and the additional first integrals
and commuting vector fields in question are also analytic), then we don't
need any Diophantine or nonresonance condition for the existence of a
convergent Poincaré-Birkhoff normalization. More precisely, we have:

\begin{thm}[\cite{Zung-Birkhoff2001,Zung-PD2002}]
\label{thm:PBmain} Let $X$ be a local analytic (non-Hamiltonian, isochore,
or Hamiltonian) vector field in $(\bbK^m,0)$ (or in $(\bbK^{2n},0)$ with a standard
symplectic structure), where $\bbK = \bbR$ or $\bbC$, such that $X(0) =
0$. Then $X$ admits a convergent Poincaré-Birkhoff normalization in a
neighborhood of $0$.
\end{thm}

Partial cases of the above theorem were obtained earlier by many
authors, including Rüssmann \cite{Russmann-NF1964} (the nondegenerate Hamiltonian
case with 2 degrees of freedom),
Vey \cite{Vey-Separable1978,Vey-Isochore1979} (the nondegenerate
Hamiltonian and isochore cases), Ito \cite{Ito-Birkhoff1989} (the
nonresonant Hamiltonian case), Ito \cite{Ito-Birkhoff1992} and
Kappeler et al. \cite{KaKoNe-Birkhoff1998} (the Hamiltonian case
with a simple resonance), Bruno and Walcher \cite{BrWa-NF1994} (the
non-Hamiltonian case with $m=2$). These
authors, except Vey who was more geometric, relied on long and heavy analytical
estimates to show the convergence of an infinite normalizing coordinate
transformation process. On the other hand, the proof of Theorem \ref{thm:PBmain}
in \cite{Zung-Birkhoff2001,Zung-PD2002}
is based on the toric point of view and is relatively short.

Following \cite{Zung-Birkhoff2001}, we will give here a sketch of the
proof of the above theorem in the Liouville-integrable case. The other
cases are similar, and of course the theorem is valid for Hamiltonian
vector fields which are integrable in generalized Liouville sense as well.
According to Proposition \ref{prop:realcomplex}, it is enough to show the
existence of a holomorphic normalization.  We will do it
 by finding local Hamiltonian ${\mathbb T}^1$-actions
which preserve the moment map of an analytically completely integrable
system. The Hamiltonian function generating such an action is an {\bf
action function}. If we find $(n-q)$ such ${\mathbb T}^1$-actions, then
they will automatically commute and give rise to a Hamiltonian ${\mathbb
T}^{n-q}$-action.

To find an action function, we will use the Mineur-Arnold formula $ P =
\int_{\Gamma} \beta$, where $P$ denotes an action function, $\beta$
denotes a primitive 1-form (i.e. $\omega = d\beta$ is the symplectic
form), and $\Gamma$ denotes an 1-cycle (closed curve) lying on a level set
of the moment map. To show the existence of such 1-cycles $\Gamma$, we
will use an approximation method, based on the existence of a formal
Birkhoff normalization.

Denote by ${\bf G} = (G_1 = H,G_2,...,G_n): ({\mathbb C}^{2n},0) \to ({\mathbb
C}^{n},0)$ the holomorphic momentum map germ of a given complex analytic
Liouville-integrable Hamiltonian system. Let $\epsilon_0 > 0$ be a small positive number such that ${\bf
G }$ is defined in the ball $\{z = (x_j,y_j) \in {\mathbb C}^{2n}, |z| <
\epsilon_0\}$. We will restrict our attention to what happens inside this ball. As
in Subsection \ref{subsection:action},
we may assume that in the symplectic coordinate system $z =
(x_j,y_j) $ we have
\begin{equation}
H = G_1 = H^s + H^n + H^{(3)} + H^{(4)} + ...
\end{equation}
with
\begin{equation}
H_s = \sum_{k=1}^{n-q} \alpha_k F^k, \, \, F^k = \sum_{j=1}^n
\rho^k_j x_jy_j ,
\end{equation}
with no resonance relations among $\alpha_1,...,\alpha_{n-q}$. We will fix this
coordinate system $z = (x_j,y_j)$, and all functions will be written in this
coordinate system.

The real and imaginary parts of the Hamiltonian vector fields of
$G_1,...,G_n$ are in involution and their infinitesimal $\bbC^n$-action
defines an {\bf associated singular Lagrangian fibration} in the ball $\{z
= (x_j,y_j) \in {\mathbb C}^{2n}, |z| < \epsilon_0\}$. For each $z$ we
will denote the fiber which contains $z$ by $M_{z}$. If $z$ is a point
such that ${\bf G} (z)$ is a regular value for the momentum map, then
$M_z$ is a connected component of ${\bf G}^{-1}({\bf G} (z))$.

Denote by
\begin{equation}
S = \{ z \in {\mathbb C}^{2n}, |z| < \epsilon_0, dG_1 \wedge dG_2 \wedge ... \wedge
dG_n (z) = 0 \}
\end{equation}
the singular locus of the moment map, which is also the set of singular
points of the associated singular foliation. What we need to know about
$S$ is that it is analytic and of codimension at least 1, though for
generic integrable systems $S$ is in fact of codimension 2. In particular,
we have the following \L ojasiewicz inequality
\cite{Lojasiewicz-Division1959}: there exist a positive number $N$ and a
positive constant $C$ such that
\begin{equation}
\label{equation:Lojasiewicz} |dG_1 \wedge ... \wedge dG_n (z) | > C (d(z, S))^N
\end{equation}
for any $z$ with $|z| < \epsilon_0$, where the norm applied to $dG_1 \wedge ...
\wedge dG_n (z)$ is some norm in the space of $n$-vectors, and $d(z,S)$ is the
distance from $z$ to $S$ with respect to the Euclidean metric. In the above
inequality, if we change the coordinate system, then only $\epsilon_0$ and $C$ have
to be changed, $N$ (the \L ojasiewicz exponent) remains the same.

We will choose an infinite decreasing series of small numbers $\epsilon_m$
($m=1,2,...$), as small as needed, with $\lim_{m \to \infty} \epsilon_m = 0$, and
define the following open subsets $U_m$ of ${\mathbb C}^{2n}$:
\begin{equation}
U_m = \{z \in {\mathbb C}^{2n}, |z| < \epsilon_m, d(z,S) > |z|^m  \}
\end{equation}

We will also choose two infinite increasing series of natural numbers $a_m$ and
$b_m$ ($m = 1,2,...)$, as large as needed, with $\lim_{m \to \infty} a_m = \lim_{m
\to \infty} b_m = \infty$. It follows from Birkhoff's formal normalization that there is a series of local holomorphic
symplectic coordinate transformations $\Phi_m$, $m \in {\mathbb N}$, such that the
following two conditions are satisfied :

a) The differential of $\Phi_m$ at $0$ is identity for each $m$, and for any two
numbers $m,m'$ with $m' > m$ we have
\begin{equation}
\Phi_{m'}(z) = \Phi_m(z) + O (|z|^{a_m}) .
\end{equation}
In particular, there is a formal limit $\Phi_{\infty} = \lim_{m \to \infty} \Phi_m$.

b) The moment map is normalized up to order $b_m$ by $\Phi_m$. More precisely, the
functions $G_j$ can be written as
\begin{equation}
G_j (z) = G_{(m)j} (z) + O (|z|^{b_m}), \, j=1,...n,
\end{equation}
with $G_{(m)j}$ such that
\begin{equation}
\{ G_{(m)j}, F_{(m)}^k \} = 0 \,\,\,\,  \forall j=1,...n, \, k=1,...,n-q .
\end{equation}
Here the functions $F_{(m)}^k$ are quadratic functions
\begin{equation}
F_{(m)}^k (x,y) = \sum_{j=1}^n \rho^k_j x_{(m)j} y_{(m)j}
\end{equation}
in local symplectic coordinates
\begin{equation}
(x_{(m)},y_{(m)}) = \Phi_m (x,y) .
\end{equation}

Notice that $F_{(m)}^k$ has the same form as $F^k$, but with respect to a
different coordinate system. When considered in the original coordinate system
$(x,y)$, $F_{(m)}^k$ is a different function than $F^k$, but the quadratic
part of $F_{(m)}^k$ is $F^k$.

Denote by $\Gamma^k_m(z)$ the orbit of the real part of the periodic Hamiltonian
vector field $X_{iF_{(m)}^k}$ which goes through $z$. Then for any $z' \in
\Gamma^k_m(z)$ we have $G_{(m)j}(z') = G_{(m)j}(z)$ and $|z'| \simeq |z|$,
therefore

\begin{equation}
\label{equation:G} |{\bf G}(z') - {\bf G}(z)| = O (|z'|^{b_m}) .
\end{equation}

(Note that we can choose the numbers $a_m$ and $b_m$ first, then choose
the radii $\epsilon_m$ of small open subsets to make them sufficiently
small with respect to $a_m$ and $b_m$, so that the equivalence $O
(|z'|^{b_m}) \simeq O (|z|^{b_m})$ makes sense). It follows from the
definition of $U_m$ and \L ojasiewicz inequalities that we also have
\begin{equation}
|d G_1 (z') \wedge ... \wedge dG_n(z')| >  d(z,S)^N > |z|^{mN}
\end{equation}
for any $z \in U_m$ and $z' \in \Gamma^k_m(z)$, provided that $b_m$ is
sufficiently large and $\epsilon_m$ is sufficiently small. Assuming that
$b_m \gg mN$, we can project the curve $\Gamma_m^k(z)$ on $M_z$ in a
unique natural way up to homotopy. Denote the image of the projection by
$\tilde{\Gamma}_m^k(z)$, and define the following function $P_m^k$ on
$U_m$:

\begin{equation}
P_m^k (z) = \oint_{\widetilde\Gamma_m^k(z)} \sum_{j=1}^n x_j dy_j\ .
\end{equation}

One then checks that $P_m^k$ is a uniformly bounded (say by 1) holomorphic
first integral of the system on $U_m$, and moreover $P_m^k$ coincides with
$P_{m'}^k$ on $U_m \cap U_{m'}$ for any $m, m'$, and hence we have a
holomorphic first integral $P^k$ on $U = \bigcup_{m=1}^\infty U_m$. The
following lemma \ref{lemma:extension} about holomorphic extension says
that $P^k$ can be extended to a holomorphic first integral of the system
in a neighborhood of 0. It is easy to see that $P^k$ is an action function
(because $P^k = \lim_{m \to \infty} \sqrt{-1} F^k_{(m)}$), i.e. its
corresponding Hamiltonian flow is periodic of period $2\pi$. Since
$k=1,\hdots,n-q$, we have $n-q$ action functions, whose flows commute and
generate the required Hamiltonian $\bbT^{n-q}$-action which preserves the
system. \QED

The following lemma on holomorphic extension, which is interesting in its
own right, implies that the action functions $P^{k}$ constructed above can
be extended holomorphically to a neighborhood of $0$.

\begin{lemma}
\label{lemma:extension} Let $U = \bigcup_{m=1}^{\infty} U_m $, with $U_m =
\{ x \in {\mathbb C}^n, |x| < \epsilon_m, d(x,S) > |x|^m \}$, where
$\epsilon_m$ is an arbitrary series of positive numbers and $S$ is a local
proper complex analytic subset of ${\mathbb C}^n$  ($codim_{\mathbb C} S
\geq 1$). Then any bounded holomorphic function on $U$ has a holomorphic
extension to a neighborhood of $0$ in ${\mathbb C}^n$.
\end{lemma}

See \cite{Zung-Birkhoff2001} for the proof of Lemma \ref{lemma:extension}.
It is a straightforward proof in the case $S$ is non-singular or is a
normal crossing, and makes use of a desingularization of $S$ in the
general case. \QED

\subsection{Torus action near a compact singular orbit} \hfill

Consider a real analytic integrable vector field $X$ on a real analytic manifold $M^m$ of
dimension $m = p+q$, with the aid of a $p$-tuple ${\bf X} =(X_1,...,X_p)$
of commuting analytic vector fields and a $q$-tuple
${\bf F} = (F_1,...,F_q)$ of analytic common first integrals:
$[X,X_i] = [X_i,X_j] = 0, X(F_j) = X_i(F_j) = 0 \ \forall i,j $.
In the Hamiltonian case, when there is an analytic Poisson structure on $M^m$, we suppose that
the system is integrable in generalized Liouville sense, i.e. the
vector fields $X, X_1,...,X_p$ are Hamiltonian.

The commuting vector fields $X_1,...,X_p$ generate an infinitesimal
$\bbR^p$-action on $M$ -- as usual, its orbits will be called orbits of the system.
The map ${\bf F}: M^m \to \bbR^q$ is constant on
the orbits of the system. Let $O \subset M^m$ be a singular orbit of
dimension $r$ of the system, $0 \leq r <
p$. We suppose that $O$ is a compact submanifold of $M^m$ (or more
precisely, of the interior of $M^m$ if $M^m$ has boundary). Then $O$ is a torus of
dimension $r$. Denote by $N$ the connected component of ${\bf F}^{-1}({\bf
F}(O))$ which contains $O$. A natural question arises: does there exist a
$\bbT^r$-action in a neighborhood of $O$ or $N$, which preserves the
system and is transitive on $O$ ?

The above question has been answered positively in
\cite{Zung-Tedemule2003}, under a mild condition called the {\it finite
type condition}. To formulate this condition, denote by $M_\bbC$ a small
open complexification of $M^m$ on which the complexification ${\bf
X}_\bbC, {\bf F}_\bbC$ of $\bf X$ and $\bf F$ exists. Denote by $N_\bbC$ a
connected component of ${\bf F}_\bbC^{-1}({\bf F}(O))$ which contains $N$.

\begin{defn} With the above notations, the singular orbit $O$ is called
of {\bf  finite type} if there is only a finite number of orbits of the
infinitesimal action of $\bbC^p$ in $N_\bbC$, and $N_\bbC$ contains a
regular point of the map $\bf F$.
\end{defn}

For example, all nondegenerate singular orbits are of finite type (see
Section \ref{section:nondegenerate}). We conjecture that every singular
orbit of an algebraically integrable system is of finite type.

\begin{thm}[\cite{Zung-Tedemule2003}]
\label{thm:tedemule-action} With the above notations, if $O$ is a compact
finite type singular orbit of dimension $r$, then there is a real analytic
torus action of $\bbT^r$ in a neighborhood of $O$ which preserves the
integrable system $({\bf X},{\bf F})$ and which is transitive on $O$. If
moreover $N$ is compact, then this torus action exists in a neighborhood
of $N$. In the Hamiltonian case this torus action also preserves the
Poisson structure.
\end{thm}

Notice that Theorem \ref{thm:tedemule-action}, together with Theorem \ref{thm:PBmain}
and the toric characterization of Poincaré-Birkhoff normalization,
provides an analytic Poincaré-Birkhoff normal form in the neighborhood a
singular invariant torus of an integrable system.

Denote by $\cA_O$ the local automorphism group of the integrable system
$({\bf X},{\bf F})$ at $O$, i.e. the group of germs of local analytic
automorphisms of $({\bf X},{\bf F})$ in vicinity of $O$ (which preserve
the Poisson structure in the Hamiltonian case). Denote by $\cA_O^0$ the
subgroup of $\cA_O$ consisting of elements of the type $g^1_Z$, where $Z$
is a analytic vector field in a neighborhood of $O$ which preserves the
system and $g^1_Z$ is the time-1 flow of $Z$. The torus in the previous
theorem is of course a Abelian subgroup of $\cA_O^0$. Actually, the
automorphism group $\cA_O$ itself is essentially Abelian in the finite
type case:

\begin{thm}[\cite{Zung-Tedemule2003}]
If $O$ is a compact finite type singular orbit as above,
then $\cA^0_O$ is an Abelian normal subgroup of
$\cA_O$, and $\cA_O/\cA^0_O$ is a finite group.
\end{thm}

The above two theorems are very closely related: their proofs are almost
the same. Let us indicate here the main
ingredients of the proof of Theorem \ref{thm:tedemule-action}:

For simplicity, we will assume that $r=1$, i.e. $O$ is a circle
(the case $r > 1$ is absolutely similar). Since $O$ is of finite
type, there is a regular complex orbit $Q$ in $N_\bbC$ of
dimension $p$ whose closure contains $O$. $Q$ is a flat affine
manifold (the affine structure is given by the $\bbC^p$-action, so
we can talk about geodesics on $Q$. If we can find a closed
geodesic $\gamma_Q$ on $Q$, then it is a periodic orbit of period
1 of a vector field of the type $\sum a_j X_j$ on $Q$ (with $a_j$
being constants) on $Q$. Since the points of $Q$ are regular for
the map $\bf F$, using implicit function theorem, we can construct
a vector field of the type $\sum a_j X_j$, with $a_j$ now being
holomorphic functions which are functionally dependent on $\bf F$
(so that this vector fields preserves the system), and which is
periodic of period 1 near $\gamma_Q$. With some luck, we will be
able to extend this vector field holomorphically to a vector field
in a neighborhood of $O$ so that $O$ becomes a periodic orbit of
it, and we are almost done: if the vector field is not
real-analytic, then its image under a complex involution will be
another periodic vector field which preserves the system; the two
vector fields commute (because the system is integrable) and we
can fabricate from them a real-analytic periodic vector field,
i.e. a real-analytic $\bbT^1$-action in a neighborhood of $O$, for
which $O$ is a periodic orbit.

The main difficulty lies in finding the closed geodesic $\gamma_Q$ (which
satisfies some additional conditions). We will do it inductively: let $
O_1 = O_\bbC \ (O \subset O_C), O_2, \dots, O_k = Q$ be a maximal chain of
complex orbits of the system in $N_\bbC$ such that $O_i$ lies in the
closure of $O_{i+1}$ and $O_i \neq O_{i+1}$. Then on each $O_i$, we will
find a closed geodesic $\gamma_i$, such that each $\gamma_{i+1}$ is
homotopic to a multiple of $\gamma_{i}$ in $O_i \cup O_{i+1}$, starting
with $\gamma_1 = O$. We will show how to go from $O = \gamma_1$ to
$\gamma_2$ (the other steps are similar). Without loss of generality, we
may assume that $O$ is a closed orbit for $X_1$. Take a small section $D$
to $O$ in $M$, and consider the Poincare map $\phi$ of $X_1$ on $D$. Let
$Y = O_2 \cap D_\bbC$. Then $Y$ is a affine manifold (whose affine
structure is projected from $O_2$ by $X_1$). Let $y$ be a point in $Y$. We
want to connect $y$ to $\phi(y)$ by a geodesic in $Y$. If we can do it,
then the sum of this geodesic segment with the orbit of $X_1$ going from
$y$ to $\phi(Y)$ can be modified into a closed geodesic $\gamma_2$ on
$O_2$. Unfortunately, in general, we cannot connect $y$ to $\phi(y)$ by a
geodesic in $Y$, because $Y$ is not ``convex''. But a lemma says that $Y$
can be cut into a finite number of convex pieces, and as a consequences
$y$ can be connected geodesically  to $\phi^N(y)$ for some power $\phi^N$
($N$-time iteration) of $\phi$. See \cite{Zung-Tedemule2003} for the
details. \QED

Theorem \ref{thm:tedemule-action} reduces the study of the behavior of
integrable systems near compact singular orbits to the study of fixed
points with a finite Abelian group of symmetry (this group arises from the
fact that the torus action is not free in general, only locally free). For
example, as was shown in \cite{Zung-Degenerate2000}, the study of {\it
corank-1} singularities of Liouville-integrable systems is reduced to the
study of families of functions on a $2$-dimensional symplectic disk which
are invariant under the rotation action of a finite cyclic group
$\bbZ/\bbZ_k$, where one can apply the theory of singularities of
functions with an Abelian symmetry developed by Wassermann
\cite{Wassermann-Symmetry1988} and other people. A (partial)
classification up to diffeomorphisms of corank-1 degenerate singularities
was obtained by Kalashnikov \cite{Kalashnikov-1Corank1998} (see also
\cite{Zung-Degenerate2000,GoSt-1Corank1987}), and symplectic invariants
were obtained by Colin de Verdière \cite{Colin-Singular2003}.

\section{Nondegenerate singularities}
\label{section:nondegenerate}

In this section, we will consider only smooth Liouville-integrable
Hamiltonian systems, though many ideas and results can probably be
extended to other kinds of integrable systems.

\subsection{Nondegenerate singular points} \hfill

Consider the momentum map
${\bf F}= (F_1,...,F_n): (M^{2n},\omega) \to \bbR^n$ of a smooth
integrable Hamiltonian system on a symplectic manifold $(M^{2n},\omega)$.
In this Section, we will forget about the original Hamiltonian function,
and study the momentum map instead.

For a point $z \in M$, denote $\rank z = \rank d{\bf F}(z)$, where $d{\bf
F}$ denotes the differential of $\bf F$. This number is equal to the
dimension of the orbit of the system (i.e. the infinitesimal Poisson
$\bbR^n$-action generated by $X_{F_1},...,X_{F_n}$) which goes through
$z$. If $\rank z < n$ then $z$ is called a {\bf singular point}. If $\rank
z = 0$ then $z$ is a {\bf fixed point} of the system.

If $z$ is a fixed point, then the quadratic parts $F^{(2)}_1,...,F^{(2)}_n$ of
the components $F_1,...,F_n$ of the momentum map at $z$ are
Poisson-commuting  and they form an Abelian subalgebra,
${A}_z$, of the Lie algebra $Q(2n,\bbR)$ of homogeneous quadratic
functions of $2n$ variables under the standard Poisson bracket.
Observe that the algebra $Q(2n,\bbR)$ is isomorphic to the
symplectic algebra $sp(2n,\bbR)$.

A fixed point $z$ will be called {\bf nondegenerate} if ${A}_z$ is a
Cartan subalgebra of $Q(2n,\mathbb R)$. In this case, according to
Williamson \cite{Williamson-1935}, there is a triple of nonnegative
integers $(k_e,k_h,k_f)$ such that $k_e + k_h + 2 k_f = n$, and a
canonical coordinate system $(x_i,y_i)$ in $\bbR^{2n}$, such that ${A}_z$
is spanned by the following quadratic functions $h_1,...,h_n$:
\begin{equation}
\label{eqn:h_i}
\begin{array}{l}
h_i = x_i^2 + y_i^2 \ \ {\rm for} \ \ 1 \leq i \leq k_e \ ;  \\
 h_i = x_iy_i \ \ {\rm for} \ \ k_e+1 \leq i \leq k_e+k_h \ ; \\
h_i = x_i y_{i+1}- x_{i+1} y_i \ \ {\rm and} \\
                            h_{i+1} = x_i y_i + x_{i+1} y_{i+1}
\ \ {\rm for} \ \ i = k_e+k_h+ 2j-1, \ 1 \leq j \leq k_f \ .
\end{array}
\end{equation}

The triple $(k_e,k_h,k_f)$ is called the {\bf Williamson type} of (the
system at) $z$. $k_e$ is the number of {\bf elliptic} components (and $
h_1,...,h_{k_e}$ are elliptic components), $k_h$ is the number of {\bf
hyperbolic} components, and $k_f$ is the number of {\bf focus-focus}
components. If $k_h = k_f = 0$ then $z$ is called an {\bf elliptic
singular point}.

The local structure of nondegenerate singular points is given by the
following theorem of Eliasson.

\begin{thm}[Eliasson \cite{Eliasson-Thesis1984,Eliasson-NF1990}]
\label{thm:Eliasson} If $z$ is a nondegenerate fixed point of a smooth
Liouville-integrable Hamiltonian system then there
is a smooth Birkhoff normalization. In other words, the singular Lagrangian
foliation given by the momentum map $\bf F$ in a neighborhood of $z$
is locally smoothly symplectomorphic to the ``linear'' singular Lagrangian fibration
given by the quadratic map $(h_1,...,h_n): \bbR^{2n} \to \bbR^n$ with the
standard symplectic structure on $\bbR^{2n}$.
\end{thm}

The elliptic case of the above theorem is also obtained independently by
Dufour and Molino \cite{DuMo-AA1991}. The case is one degree of freedom is
due to Colin de Verdière and Vey \cite{CoVe-Isochore1979}. The analytic
case of the above theorem is due to Vey \cite{Vey-Separable1978}, and is
superseded by Theorem \ref{thm:PBmain}. There is also a semiclassical
version of Eliasson's theorem (quantum Birkhoff normal form), which is due
to Vu Ngoc San \cite{San-Birkhoff2000}.

The proof of Eliasson's theorem \cite{Eliasson-Thesis1984,Eliasson-NF1990}
is quite long and highly technical: The first step is to use division
lemmas in singularity theory to show that the local singular fibration
given by the momentum map is diffeomorphic (without the symplectic
structure) to the linear model. Then one uses a combination of averaging,
Moser's path method, and technics similar to the ones used in the proof of
Sternberg's smooth linearization theorem for vector fields, to show that
the symplectic form can also be normalized smoothly. In fact, Eliasson's
proof of his theorem is not quite complete, except for the elliptic case,
because it lacks some details which were difficult to work out, see
\cite{MiSa-Singular2004}.

A direct consequence of Eliasson's theorem is that, near a nondegenerate
fixed point of Williamson type $(k_e,k_h,k_f)$, there is a local smooth
Hamiltonian $\bbT^{k_e+k_f}$-action which preserves the system: each
elliptic or focus-focus component provides one $\bbT^1$-action. In the
analytic case, Birkhoff normalization gives us a $\bbT^{n}$-action, but it
acts in the complex space, and in the real space we only see a
$\bbT^{k_e+k_f}$-action.

\subsection{Nondegenerate singular orbits} \hfill

Let $x \in M$ be a singular point of $\rank x = m \geq 0$. We may assume
without loss of generality that $dF_1 \wedge ... \wedge dF_m (x) \neq 0$,
and a local symplectic reduction near $x$ with respect to the local free
$\bbR^m$-action generated by the Hamiltonian vector fields
$X_{F_1},...,X_{F_m}$ will give us an $m$-dimensional family of local
integrable Hamiltonian systems with $n-m$ degrees of freedom. Under this
reduction, $x$ will be mapped to a fixed point in the reduced system, and
if this fixed point is nondegenerate according to the above definition,
then $x$ is called a {\bf nondegenerate singular point} of {\bf rank} $m$
and {\bf corank} $(n-m)$. In this case, we can speak about the Williamson
type $(k_e,k_h,k_f)$ of $x$, and we have $k_e + k_h + 2 k_f = m$.

A {\bf nondegenerate singular orbit} of the system is an orbit (of the
infinitesimal Poisson $\bbR^n$-action) which goes through a nondegenerate
singular point. Since all points on a singular orbit have the same
Williamson type, we can speak about the Williamson type and the corank of
a nondegenerate singular orbit. We have the following generalization of
Theorem \ref{thm:Eliasson} to the case of compact nondegenerate singular
orbits:

\begin{thm}[Miranda--Zung \cite{MiZu-Orbit2003}]
\label{thm:MZ-Orbit} If $O$ is a compact nondegenerate singular orbit of a
smooth Liouville-integrable Hamiltonian system, then the singular
Lagrangian fibration given by the momentum map in a neighborhood of $O$ is
smoothly symplectomorphic to a linear model. Moreover, if the system is
invariant under a symplectic action of a compact Lie group $G$ in a
neighborhood of $O$, then the above smooth symplectomorphism to the linear
model can be chosen to be $G$-equivariant.
\end{thm}

The linear model in the above theorem can be constructed as follows:
Denote by $(p_1,\hdots,p_m)$ a linear coordinate system of a small ball
$D^m$ of dimension $m$, $(q_1 (mod\ 1),\hdots,q_m (mod\ 1))$ a standard
periodic coordinate system of the torus $\bbT^m$, and
$(x_1,y_1,\hdots,x_{n-m},y_{n-m})$ a linear coordinate system of a small
ball $D^{2(n-m)}$ of dimension $2(n-m)$. Consider the manifold
\begin{equation}
\label{eqn:V} V = D^m \times \bbT^m \times D^{2(n-m)}
\end{equation}
 with the standard symplectic form $\sum dp_i
\wedge dq_i + \sum dx_j \wedge dy_j$, and the following momentum map: $
\label{eqn:linearmm} ({\bf p},{\bf h}) =
(p_1,\hdots,p_m,h_1,\hdots,h_{n-m}): V \rightarrow \bbR^n ,$ where
$(h_1,\hdots,h_{n-m})$ are quadratic functions given by Equation
($\ref{eqn:h_i}$). A symplectic group action on $V$ which preserves the
above momentum map is called {\it linear} if it on the product $V = D^m
\times \bbT^m \times D^{2(n-m)}$ componentwise, the action  on $D^m$ is
trivial, the action on $\bbT^m$ is by translations with respect to the
coordinate system $(q_1,\hdots,q_m)$, and the action on $D^{2(n-m)}$ is
linear.

Let $\Gamma$ be a finite group with a free linear symplectic
action $\rho(\Gamma)$ on $V$ which preserves the momentum map.
Then we can form the quotient integrable system with the momentum map
\begin{equation}
\label{eqn:twistedmodel} ({\bf p},{\bf h}) =
(p_1,\hdots,p_m,h_1,\hdots,h_{n-m}): V/\Gamma \rightarrow \bbR^n \ .
\end{equation}
The set $\{p_i=x_i=y_i = 0\} \subset V/\Gamma$ is a compact orbit of
Williamson type $(k_e,k_f,k_h)$ of the above system. The above system on
$V/\Gamma$ is called the {\bf linear model} of Williamson type
$(k_e,k_f,k_h)$ and twisting group $\Gamma$, or more precisely, twisting
action $\rho(\Gamma)$. (It is called a direct model if $\Gamma$ is
trivial, and a twisted model if $\Gamma$ is nontrivial). A symplectic
action of a compact group $G$ on $V/\Gamma$ which preserves the momentum
map $(p_1,\hdots,p_m,h_1,\hdots,h_{n-m})$ is called linear if it comes
from a linear symplectic action of $G$ on $V$ which commutes with the
action of $\Gamma$.

The case with $G$ trivial and $n=2, k_h=1, k_e=k_f=0$  of Theorem \ref{thm:MZ-Orbit}
is due to Colin de Verdière and Vu Ngoc San \cite{CoVu-2D2002},
and independently Curr\'as-Bosch and Miranda \cite{CuMi-Linearization2002}.
A direct consequence of Theorem \ref{thm:MZ-Orbit} is that the group of
local smooth symplectic automorphisms of a smooth Liouville-integrable
system near a compact nondegenerate singular orbit is Abelian, see
\cite{MiZu-Orbit2003}.

\subsection{Nondegenerate singular fibers} \hfill

In this subsection, we will assume that the momentum map ${\bf F}: M^{2n}
\to \bbR^n$ is proper. A singular connected component of a level set of
the momentum map will be called a {\bf singular fiber} of the system. A
singular fiber may contain one orbit (e.g. in the elliptic nondegenerate
case), or many orbits, some of them singular and some of them regular. A
singular fiber $N_c$ is called {\bf nondegenerate} if any point $z \in
N_c$ is either regular or nondegenerate singular. The nondegeneracy is an
open condition: if a singular fiber is nondegenerate then nearby singular
fibers are also nondegenerate.

By a {\bf singularity} of a Liouville-integrable system, we mean the germ
of the system near a singular fiber, together with the  symplectic form
and the Lagrangian fibration. We will denote a singularity by
$(\cU(N_c),\omega,\cL)$, where $\cU(N_c)$ denotes a small ``tubular''
neighborhood of $N_c$, and $\cL$ denotes the Lagrangian fibration. If
$N_c$ is nondegenerate then $(\cU(N_c),\omega,\cL)$ is also called
nondegenerate.

A simple lemma \cite{Zung-AL1996} says that if $N_c$ is a nondegenerate
singular fiber, then all singular points of maximal corank in $N_c$ have
the same Williamson type. We define the rank and the Williamson
type of a nondegenerate singularity $(\cU(N_c),\omega,\cL)$ to be the
rank and the Williamson type of a singular point of maximal corank in
$N_c$.

The following theorem may be viewed as the generalization of
Liouville--Mineur--Arnold theorem to the case of nondegenerate singular
fibers:

\begin{thm}[\cite{Zung-AL1996}]
\label{thm:nondegenerate}
Let $(\cU(N_c),\omega,\cL)$ be a nondegenerate smooth singularity of rank $m$
and Williamson type $(k_e,k_h,k_f)$ of a Liouville-integrable system with a proper
momentum map. Then we have: \\
a) There is effective Hamiltonian $\bbT^{m+k_e+k_f}$-action in
$(\cU(N_c),\omega,\cL)$ which preserves the system. The dimension $m+k_e+k_f$
is maximal possible. There is a locally free $\bbT^m$-subaction of this action. \\
b) There is a partial action-angle coordinate system. \\
c) Under a mild additional condition,  $(\cU(N_c), \cL)$ is topologically
equivalent to an almost direct product of simplest (corank 1 elliptic or
hyperbolic and corank 2 focus-focus) singularities.
\end{thm}

Assertion b) of the above theorem means that we can write
$(\cU(N_c),\omega)$ as $(D^{m} \times \bbT^m
\times P^{2k})  / \Gamma$ with
\begin{equation}
\omega = \sum_1^m dp_i \wedge dq_i + \omega_1
\end{equation}
where $\omega_1$ is a symplectic form on $P^{2k}$, the finite group
$\Gamma$ acts on the product component-wise, its action is linear on $\bbT^m$,
and the momentum map $\bf F$ does not depend on the variables $q_1,...,q_m$.

The additional condition in Assertion c) prohibits the bifurcation
diagram (i.e. the set of singular values of the momentum map) from
having ``pathologies'', see \cite{Zung-AL1996}, and it's satisfied
for all nondegenerate singularities of physical integrable systems
met in practice. The almost direct product means a product of the
type
\begin{equation}
\label{eqn:almostdirectproduct}
({\mathcal T}^{2m} \times {\mathcal E}^2_1 \times ... \times {\mathcal E}^2_{k_e}
\times  {\mathcal H}^2_1 \times ... \times {\mathcal H}^2_{k_h} \times
{\mathcal F}^4_1 \times ... \times {\mathcal F}^4_{k_f}) / \Gamma
\end{equation}
where ${\mathcal T}^{2m}$ is the germ of $(D^m \times \bbT^m, \sum_1^m dp_i \wedge
dq_i)$ with the standard Lagrangian torus fibration;
${\mathcal E}^2_i, {\mathcal F}^2_i$ and ${\mathcal H}^4_i$  are elliptic,
hyperbolic and focus-focus singularities of integrable systems
on symplectic manifolds of dimension 2, 2 and 4 respectively;
the finite group $\Gamma$  acts freely and component-wise. Remark that, in
general, a nondegenerate singularity is only topologically equivalent, but
not symplectically equivalent, to an almost direct product singularity.

The above almost direct product may remind one of the decomposition of
algebraic reductive groups into almost direct products of simple groups
and tori: though the two objects are completely different, there are some
common ideas behind them, namely infinitesimal direct decomposition, and
twisting by a finite group.

\subsection{Focus-focus singularities} \hfill

The singularities  ${\mathcal E}^2_i, {\mathcal H}^2_i, {\mathcal F}^4_i$
in (\ref{eqn:almostdirectproduct}) may be called {\bf elementary}
nondegenerate singularities; they are characterized by the fact that $k_e
+ k_h + k_f = 1$ and $\rank = 0$. Among them, elliptic singularities
${\mathcal E}^2_i$ are very simple: each elementary elliptic singularity
is isomorphic to a standard linear model (a harmonic oscillator).
Elementary hyperbolic singularities ${\mathcal H}^2_i$ are also relatively
simple because they are given by hyperbolic singular level sets of Morse
functions on 2-dimensional symplectic surfaces. On the other hand, {\it
focus-focus} singularities ${\mathcal F}^4_i$  live in 4-dimensional
symplectic manifolds, so their topological structure is somewhat more
interesting. Let us mention here some results about the structure of these
4-dimensional focus-focus singularities, see
\cite{Zung-AL1996,Zung-FocusII2002} and references therein for more
details.

One of the most important facts about focus-focus singularities is the
existence of a $\bbT^1$-action (this is a special case of Assertion a) of
Theorem \ref{thm:nondegenerate}); many other important properties are
consequences of this $\bbT^1$-action. In fact, in many integrable systems
with a focus-focus singularity, e.g. the spherical pendulum and the
Lagrangian top, this $\bbT^1$-action is the obvious rotational symmetry,
though in some systems, e.g. the Manakov integrable system on $so(4)$,
this local $\bbT^1$-action is ``hidden''. Dynamically speaking, a
focus-focus point is roughly an unstable equilibrium point with a
$\bbT^1$- symmetry.

Each focus-focus singularity has only one singular fiber: the focus- focus
fiber, which is homeomorphic to a {\bf pinched torus} (take a torus, and
$\ell$ parallel homotopically non-trivial simple closed curves on it,
$\ell \geq 1$, then collapse each of these curves into one point). This
fact was known to Lerman and Umanskij \cite{LeUm-2D1987,LeUm-2D1988}.

From the topological point of view, we have a singular torus fibration in
a four-dimensional manifold with one singular fiber. These torus
fibrations have been studied by Matsumoto and other people, see e.g.
\cite{Matsumoto-Torus1989} and references therein, and of course the case
with a singular fiber of focus-focus type is included in their topological
classification. In particular, the number of pinches $\ell$ is the only
topological invariant. The monodromy of the torus fibration (over a
punched 2-dimensional disk) around the focus-focus fiber is given by the
matrix $\begin{pmatrix}  1 & \ell \\ 0 & 1 \end{pmatrix}$. By the way, the
case with $\ell > 1$ is topologically an $\ell$-sheet covering of the case
with $\ell=1$, and a concrete example with $\ell=1$ is the unstable
equilibrium of the usual spherical pendulum.  This phenomenon of {\it
nontrivial monodromy} (of the foliation by Liouville tori) was first
observed by Duistermaat and Cushman \cite{Duistermaat-AA1980}, and then by
some other people for various concrete integrable systems.
Now we have many different ways to look at this monodromy: from the purely
topological point of view (using Matsumoto's theory
\cite{Matsumoto-Torus1989}), from the point of view of Picard-Lefschetz
theory (see Audin \cite{Audin-Monodromy2002} and references therein), or
as a consequence of Duistermaat-Heckman formula with respect to the
above-mentioned $\bbT^1$-action (see \cite{Zung-FocusII2002}).
Quantization of focus-focus singularities leads to {\it quantum}
monodromy, see Vu Ngoc San \cite{San-Focus2000} and Section
\ref{section:global}.

Similar results, including the existence of a $\bbT^1$-action, for
focus-focus singularities of {\it non-Hamiltonian} integrable systems,
have been obtained by Cushman an Duistermaat \cite{CuDu-Focus2001}, see
also \cite{Zung-FocusII2002}.

\section{Global aspects of local torus actions}
\label{section:global}

\subsection{Sheaf of local $\bbT^1$-actions} \hfill
\label{subsection:monodromy}

Consider a smooth proper integrable system on a manifold $M$ with a given
$p$-tuple of commuting vector fields ${\bf X} = (X_1,...,X_p)$ and
$q$-tuple of common first integrals ${\bf F} = (F_1,...,F_q)$.

We will call the space of connected components of the level sets of the
map $\bf F$ the {\bf base space} of the integrable system, and denote it
by $\cB$. Since the system is proper, the space $\cB$ with the induced
topology from $M$ is a Hausdorff space. We will denote by ${\bf P}: M \to
\cB$ the projection map from $M$ to $\cB$.

For each open set $U$ of $\cB$, denote by $\cR(U)$ the set of all
$\bbT^1$-actions in ${\bf P}^{-1}(U)$ which preserve the integrable system
$({\bf X},{\bf F})$ (in the Hamiltonian case, due to generalized
Liouville-Mineur-Arnold theorem, elements of $\cR(U)$ will automatically
preserve the Poisson structure). $\cR(U)$ is an Abelian group: if two
elements $\rho_1,\rho_2$ of $\cR(U)$ are generated by two periodic vector
fields $Y_1,Y_2$ respectively, then $Y_1$ will automatically commute with
$Y_2$, and the sum $Y_1+Y_2$ generates another $\bbT^1$-action which can
be called the sum of $\rho_1$ and $\rho_2$. Actually, $\cR(U)$ is a free
Abelian group, and its dimension can vary from $0$ to $p$ (the dimension
of a regular invariant torus of the system), depending on $U$ and on the
system. If $U$ is a small disk in the regular region of $\cB$ then
$\dim_\bbZ \cR(U) = p$.

The association $U \mapsto \cR(U)$ forms a free Abelian sheaf $\cR$ over
$\cB$, which we will call the {\bf toric monodromy sheaf} of the system.
This sheaf was first introduced in \cite{Zung-IntegrableII2001} for the
case of Liouville-integrable systems, but its generalization to the cases
of non-Hamiltonian integrable systems and integrable systems in
generalized Liouville sense is obvious.

If we restrict $\cR$ to the regular region $\cB_0$ of $\cB$ (the set of
regular invariant tori of the system), then $\cB$ is a locally trivial
free Abelian sheaf of dimension $p$ (one may view it as a $\bbZ^p$-bundle
over $\cB_0$), and its monodromy (which is a homomorphism from the
fundamental group $\pi_1(\cB_0)$ of $\cB_0$ to $GL(p,\bbZ)$) is nothing
but the topological monodromy of the torus fibration of the regular part
of the system. This topological monodromy, in the case of
Liouville-integrable system, is known as the monodromy in the sense of
Duistermaat \cite{Duistermaat-AA1980}, and it is a topological obstruction
to the existence of global action-angle variables. In the case of
Liouville-integrable systems with only nondegenerate elliptic
singularities, studied by Boucetta and Molino \cite{BoMo-Elliptic1989},
$\cR$ is still a locally free Abelian sheaf of dimension $p = {1 \over 2}
\dim M$.

When the system has non-elliptic singularities, the structure of $\cR$ can
be quite complicated, even locally, and it contains a lot more information
than the monodromy in the sense Duistermaat. For example, in the case of
2-degree-of-freedom Liouville-integrable systems restricted to isoenergy
3-manifolds, $\cR$ contains information on the ``marks'' of the so-called
{\bf Fomenko-Zieschang invariant}, which is a complete topological
invariant for such systems, see e.g.
\cite{FoZi-Invariant1991,BoFo-Integrable1999}. In fact, as found out by
Fomenko, these isoenergy 3-manifolds are graph-manifolds, so the classical
theory of graph-manifolds can be applied to the topological study of these
2-degree-of-freedom Liouville-integrable systems. A simple explanation of
the fact that these manifolds are graph-manifolds is that they admit local
$\bbT^1$-actions.

The second cohomology group $H^2(\cB,\cR)$ plays an important role in the
global topological study of integrable systems, at least in the
Liouville-integrable case. In fact, if two Liouville-integrable
Hamiltonian systems have the same base space, the same singularities, and
the same toric monodromy sheaf, then their remaining topological
difference can be characterized by an element in $H^2(\cB,\cR)$, called
the (relative) {\bf Chern class}. We refer to \cite{Zung-IntegrableII2001}
for a precise definition of this Chern class for Liouville-integrable
Hamiltonian systems (the definition is quite technical when the system has
non-elliptic singularities), and the corresponding topological
classification theorem. In the case of systems without singularities or
with only elliptic singularities, this Chern class was first defined and
studied by Duistermaat \cite{Duistermaat-AA1980}, and then by
Dazord--Delzant \cite{DaDe-AA1987} and Boucetta--Molino
\cite{BoMo-Elliptic1989}.

\subsection{Affine base space, integrable surgery, and convexity}\hfill

In the case of Liouville-integrable systems, the base space $\cB$ has a
natural stratified integral affine structure (local action functions of
the system project to integral affine functions on $\cB$), and the
structure of the toric monodromy sheaf $\cR$ can be read off the affine
structure of $\cB$, see \cite{Zung-IntegrableII2001}.

The integral affine structure on $\cB$ plays an important role in the
problem of quantization of Liouville-integrable systems. A general idea,
supported by recent works on quantization of integrable systems, see e.g.
\cite{San-Focus2000,CoVu-2D2002,NeSaZi-Fractional2002}, is that one can
think of Bohr-Sommerfeld or quasi-classical quantization as a
discretization of the integral affine structure of $\cB$: after
quantization, in place of a stratified integral affine manifold, we get a
``stratified nonlinear lattice'' (of joint spectrum of the system). The
monodromy of this joint spectrum stratified lattice of the quantized
system (called {\it quantum monodromy}) naturally resembles the monodromy
of the classical system.

The idea of {\bf integrable surgery}, introduced in
\cite{Zung-IntegrableII2001}, is as follows: if we look at integrable
systems from differential topology point of view (singular torus
fibrations), instead of dynamical point of view (quasi-periodic flows),
then we can perform surgery on them in order to study their properties and
obtain new integrable systems from old ones. The surgery is first
performed at the base space level, and then lifted to the phase space. As
a side product, we also obtain new symplectic manifolds from old ones.

Let us indicate here a few interesting results obtained in
\cite{Zung-IntegrableII2001,Zung-FocusII2002} with the help of integrable
surgery:

\begin{itemize}
\item A very simple example of an exotic symplectic space $\bbR^{2n}$
(which is diffeomorphic to a standard symplectic space $(\bbR^{2n},\sum
dx_i \wedge dy_i)$ but cannot be symplectically embedded into
$(\bbR^{2n},\sum dx_i \wedge dy_i)$).

\item Construction of integrable systems on symplectic manifolds
diffeomorphic to K3 surfaces, and on other symplectic 4-manifolds.

\item Existence of a fake base space, i.e. a stratified integral affine
manifold which can be realized locally as a base space of a
Liouville-integrable system but globally cannot.

\item A new simple proof \cite{Zung-FocusII2002} of the monodromy formula
around focus-focus singularities (and their degenerate analogs) using
Duistermaat--Heckman formula. The idea itself goes beyond focus-focus
singularities and can be applied to other situations as well, for example
in order to obtain fractional monodromy \cite{NeSaZi-Fractional2002} via
Duistermaat-Heckman formula.
\end{itemize}

Integrable surgery was recently adopted by Symington in her work on
symplectic 4-manifolds \cite{Symington-Blowdown2001,Symington-4dim2002},
and in the work of Leung and Symington \cite{LeSy-AlmostToric2003} where
they gave a complete list of diffeomorphism types of 4-dimensional closed
``almost toric'' symplectic manifolds. It seems that integrable surgery
was also used by Kontsevich and Soibelman in a recent work on mirror
symmetry \cite{KoSo-Mirror2004}.

For a noncommutatively integrable Hamiltonian system on symplectic
manifolds, local action functions defined by the Mineur-Arnold formula
still project to local functions on the base space $\cB$, but since the
number of independent action functions is equal to the dimension of
invariant tori and is smaller than the dimension of $\cB$, they don't
define an affine structure on $\cB$, but rather an integral affine
structure transverse to a singular foliation in $\cB$, and under some
properness condition this transverse affine structure projects to an
affine structure on a quotient space $\widehat{\cB}$ of $\cB$, which has
the same dimension as that of invariant tori and which may be called the
{\bf reduced affine base space}.

A particular situation where $\widehat{\cB}$ looks nice is the case of
noncommutatively integrable systems generated by Hamiltonian actions of
compact Lie groups, or more generally of {\it proper symplectic
groupoids}, see \cite{Zung-Proper2004}. It was shown in
\cite{Zung-Proper2004} that in this case $\widehat{\cB}$ is an integral
affine manifold with locally convex polyhedral boundary, and we have a
kind of {\bf intrinsic convexity} from which one can recover various
convexity theorems for momentum maps in symplectic geometry, including,
for example:

\begin{itemize}
\item Atiyah--Guillemin--Sternberg--Kirwan convexity theorem
\cite{Atiyah-Convexity1982,GuSt-Convexity1982,Kirwan-Convexity1984} which
says that if $G$ is a connected compact Lie group which acts Hamiltonianly
on a connected compact symplectic manifold $M$ with an equivariant
momentum map $\mu: M \rightarrow \fg^*$ then $\mu(M) \cap \ft^*_+$ is a
convex polytope, where $\ft^*_+$ denotes a Weyl chamber in the dual of a
Cartan subalgebra of the Lie algebra of $G$.

\item Flaschka--Ratiu's convexity theorem for momentum maps of Poisson
actions of compact Poisson-Lie groups \cite{FlRa-Convexity1996}.

\item If one works in an even more general setting of Hamiltonian spaces
of proper quasi-symplectic groupoids \cite{Xu-Momentum2003}, then one also
recovers Alekseev--Malkin--Meinrenken's convexity theorem for group-valued
momentum maps \cite{AMM-GroupMoment1998}.
\end{itemize}

In order to further incite the reader to read \cite{Zung-Proper2004}, let
us mention here a recent beautiful convexity theorem of Weinstein
\cite{Weinstein-Noncompact2001} which has a clear meaning in Hamiltonian
dynamics and which fits well in the above framework of proper symplectic
groupoid actions and noncommutatively integrable Hamiltonian systems:

\begin{thm}[Weinstein \cite{Weinstein-Noncompact2001}]
\label{thm:Weinstein_convexity} For any positive-definite quadratic
Hamiltonian function $H$ on the standard symplectic space $\bbR^{2k}$,
denote by $\phi(H)$ the $k$-tuple $\lambda_1 \leq \hdots \leq \lambda_k$
of frequencies of $H$ ordered non-decreasingly, i.e. $H$ can be written as
$H = \sum \lambda_i(x_i^2 + y_i^2)$ in a canonical coordinate system. Then
for any two given positive nondecreasing $n$-tuples $\lambda =
(\lambda_1,\hdots,\lambda_k)$ and $\gamma = (\gamma_1,\hdots,\gamma_k)$,
the set
\begin{equation}
\Phi_{\lambda,\gamma} = \{\phi(H_1+H_2) \ | \ \phi(H_1) = \lambda,
\phi(H_2) = \gamma \}
\end{equation}
is a closed, convex, locally polyhedral subset of $\bbR^k$.
\end{thm}

Note the above set $\Phi_{\lambda,\gamma}$ is closed but not bounded. For
example, when $k =1$ then $\Phi_{\lambda,\gamma}$ is a half-line.

\subsection{Localization formulas}  \hfill

A general idea in analysis and geometry is to express global invariants
in terms of local invariants, via {\it localization formulas}.

Various global topological invariants, including the Chern classes (of the
tangent bundle), of the symplectic ambient manifold of a
Liouville-integrable system, can be localized at singularities of the
system. Some results in this direction can be found in recent papers of
Gross \cite{Gross-TopMirror2001} and Smith \cite{Smith-Torus2001}, though
much still waits to be worked out for general integrable systems. For
example, consider a 4-dimensional symplectic manifold with a proper
integrable system whose fixed points are nondegenerate. Then to find $c_2$
(the Euler class) of the manifold, one simply needs to count the number of
fixed points with signs: the plus sign for elliptic-elliptic ($k_e = 2,
k_h=k_f=0$ in Williamson type), hyperbolic-hyperbolic ($k_h = 2)$ and
focus-focus points, and the minus sign for elliptic-hyperbolic
($k_e=k_h=1$) points.

In symplectic geometry, there is a famous localization formula for
Hamiltonian torus actions, due to Duistermaat and Heckman
\cite{DuHe-Localization1982}. There is a topological version of this
formula, in terms of equivariant cohomology, due to Atiyah--Bott
\cite{AtBo-Equivariant1984} and Berline--Vergne
\cite{BeVe-Equivariant1982}, and a non-Abelian version due to Witten
\cite{Witten-Gauge1992} and Jeffrey--Kirwan \cite{JeKi-Nonabelian1995}. We
refer to
\cite{Audin-Torus1991,Duistermaat-Equivariant1994,GuGiKa-Moment2002} for
an introduction to these formulas. It would be nice to have analogs of
these formulas for proper integrable systems.

\section{Infinite-dimensional torus actions}

A general idea is that proper infinite dimensional integrable systems
admit infinite-dimensional torus actions. Consider, for example, the KdV
equation
\begin{equation}
u_t = -u_{xxx} + 6uu_x
\end{equation}
with periodic boundary condition $u(t,x+1) = u(t,x)$. We will view this
KdV equation as a flow on the space of functions $u$ on $\bbS^1 =
\bbR/\bbZ$. Then it is a Hamiltonian equation with the Hamiltonian
\begin{equation}
H = \int_{\bbS^1} ({1 \over 2}u_x^2 + u^3) dx
\end{equation}
and the Poisson structure ${d \over dx}$, see e.g. \cite{KaPo-KdVKAM2003}.
In other words, if $F$ and $G$ are two functional on the space of
functions on $\bbS^1$ then their Poisson bracket is
\begin{equation}
\{F,G\} = \int_{\bbS^1} {\partial F \over \partial u} {d \over dx}
{\partial G \over \partial u} dx \ .
\end{equation}
The Poisson structure ${d \over dx}$ admits a Casimir function
\begin{equation}
[u] = \int_{\bbS^1}u dx \ .
\end{equation}
The Sobolev space
\begin{equation}
\cH^1_0 = \left\{ u: \bbS^1 \rightarrow \bbR^1 \ | \int_{\bbS^1} (u_x^2 +
u^2) dx < \infty, \ [u] = 0 \right\}
\end{equation}
together with the Poisson structure $d \over dx$ is a \emph{weak}
symplectic infinite-dimensional manifold, which is symplectomorphic, via
the Fourier transform, to the symplectic Hilbert space
\begin{equation}
{\frak h}_{3/2} := \left\{(x_n,y_n)_{n \in \bbN} \ | \ x_n,y_n \in \bbR,
\sum n^3 x_n^2 + \sum n^3y_n^2  < \infty \right\}
\end{equation}
with the \emph{weak} symplectic structure $\omega = \sum dx_n \wedge
dy_n$. We have the following Birkhoff normal form theorem for the periodic
KdV equation, due to Kappeler and his collaborators Bättig, Bloch,
Guillot, Mityagin, Makarov, see \cite{BBGK-NF1995,KaPo-KdVKAM2003} and
references therein:

\begin{thm}[Kappeler et al.]
\label{thm:KdV} There is a bi-analytic 1-1 symplectomorphism $\Psi: {\frak
h}_{3/2} \rightarrow \cH^1_0  $, such that the coordinates $(x,y)$ become
global Birkhoff coordinates for the KdV equation under this
symplectomorphism, i.e. the transformed Hamiltonian $H_\Psi = H \circ
\Psi$ depend only on $x_n^2 + y_n^2, n \in \bbN$.
\end{thm}

See \cite{KaPo-KdVKAM2003} for a more precise and general statement of the
above theorem. A direct consequence of Theorem \ref{thm:KdV} is that we
have a Hamiltonian infinite-dimensional torus action generated by the
action functions $I_n = x_n^2 + y_n^2$ which preserve the KdV system and
whose orbits are exactly the level sets of the KdV. In fact, these action
functions are also found by the Mineur-Arnold integral formula.

Remark that the level sets $N_c := \{I_n = c_n\ \forall n \in \bbN\}$,
where $c_n$ are nonnegative constants such that $\sum n^3c_n < \infty$,
are {\it compact} with respect to the induced norm topology. They are
\emph{not} submanifolds of the phase space, but can still be viewed as
(infinite-dimensional) Liouville tori. The natural topology that we have
to put on the infinite-dimensional torus is the \emph{product topology},
and then it becomes a compact topological group by Tikhonoff theorem, and
the Hamiltonian action $\bbT^\infty \times \cH^1_0 \rightarrow \cH^1_0$ is
a continuous action.

So basically the periodic KdV is just an infinite-dimensional oscillator.
The periodic defocusing NLS (non-linear Schrödinger) equation is similar
and is also an infinite-dimensional oscillator, see
\cite{BBGK-NF1995,GrKaPo-NLS2003}.

There are other integrable equations, like the sine-Gordon equation and
the focusing NLS equation, which are topologically very different from
infinite-dimensional oscillators: they admit \emph{unstable}
singularities. Their topological structure was studied to some extent by
McKean, Ercolani, Forest, McLaughlin and many other people, see e.g.
\cite{McKean-Sine1981,ErMl-Topology1991,LiMl-NLS1994} and references
therein.

A general idea is that, even when an integrable PDE admits unstable
phenomena, locally near each level set the system can be decomposed into 2
parts: an unstable part which is finite-dimensional, and an
infinite-dimensional oscillator. The reason is that the energy is finite,
and to get something unstable one needs big energy, so one can only get
finitely many unstable things, the rest is just a small
(infinite-dimensional) oscillator. This idea reduces the topological study
of infinite-dimensional integrable systems to that of finite-dimensional
systems.

A concrete case study, namely the symplectic topology of the focusing NLS
equation,
\begin{equation}
-iq_t = q_{xx} + 2 \bar{q}q^2\ ,
\end{equation}
where $q$ is a complex-valued function on $\bbS^1$ for each $t$, is the
subject of a joint work in progress of Thomas Kappeler, Peter Topalov and
myself \cite{KaToZu-NLS2004}. It is a Hamiltonian system on the Sobolev
space $\cH^1(\bbS^1,\bbC)$ of complex-valued functions on the circle with
the Poisson bracket
\begin{equation}
\{F,G\} = i \int_{\bbS^1} \left({\partial F \over \partial q} {\partial G
\over
\partial \bar{q}} - {\partial G \over \partial q} {\partial F \over
\partial \bar{q}}\right) dx \ ,
\end{equation}
and its Hamiltonian is
\begin{equation}
H = \int_{\bbS^1} (q_x \bar{q}_x - q^2\bar{q}^2) dx \ .
\end{equation}

Here is our speculation as to what happens there: For any $q$ which
belongs to a dense open subset ${\mathcal M}_0$ of $\cH^1(\bbS^1,\bbC)$
(the ``almost-regular set''), the level set $Iso(q)$ (which is the same as
the isospectral set of the Zakharov-Shabat operator) is a torus (of
infinite dimension in general), and there is a neighborhood $\cU(Iso(q))$
of $Iso(q)$ in the phase space which admits a full action-angle system of
coordinates similar to the KdV case. For any $q \notin {\mathcal M}_0$, a
neighborhood $\cU(Iso(q))$ of $Iso(q)$ still admits a torus action of
``finite corank'' and we have a partial Birkhoff coordinate system of
finite codimension. If $q \notin {\mathcal M}_0$ is ``nondegenerate'' (it
corresponds to a condition on the spectrum of $q$), then $\cU(Iso(q))$
together with the fibration by the level sets can be written topologically
as a direct product of an infinite-dimensional oscillator and a finite
number of focus-focus singularities. The fact that unstable nondegenerate
singularities of the focusing NLS are of focus-focus type can be seen from
the work of Li and McLaughlin \cite{LiMl-NLS1994}, and is probably due to
the $\bbT^1$-symmetry of the system (translations in $x$-variable). If we
restrict the focusing NLS system to even functions $(q(-x) = q(x))$, then
it is still an integrable Hamiltonian system, but with hyperbolic instead
of focus-focus singularities.


\providecommand{\bysame}{\leavevmode\hbox to3em{\hrulefill}\thinspace}

\end{document}